\documentclass[aoas]{imsart}

\RequirePackage{amsthm,amsmath,amsfonts,amssymb}
\RequirePackage[authoryear]{natbib}

\startlocaldefs
\theoremstyle{plain}
\newtheorem{theorem}{Theorem}
\theoremstyle{remark}
\newtheorem{assumption}{Assumption}
\newtheorem{lemma}{Lemma}
\newtheorem{remark}{Remark}
\usepackage{url}
\usepackage{graphicx}
\usepackage{prodint}
\usepackage[shortlabels]{enumitem} 
\usepackage{booktabs}
\usepackage{tikz}
\usetikzlibrary{arrows,positioning,calc, external} 
\usepackage{float}

\makeatletter
\newcommand*\bigcdot{\mathpalette\bigcdot@{.5}}
\newcommand*\bigcdot@[2]{\mathbin{\vcenter{\hbox{\scalebox{#2}{$\m@th#1\bullet$}}}}}
\makeatother

\DeclareMathOperator*{\argmax}{arg\,max}
\newcommand*\diff{\mathop{}\!\mathrm{d}}

\tikzset{
    >=stealth',
    punkt/.style={
           rectangle,
           rounded corners,
           draw=black, thick,
           text width=7em,
           minimum height=2em,
           text centered},
    punkts/.style={
           rectangle,
           rounded corners,
           draw=black, thick,
           text width=3em,
           minimum height=2em,
           text centered},       
    punktl/.style={
           re
           tangle,
           rounded corners,
           draw=black, thick,
           
           text width=7em,
           minimum height=2em,
           text centered},
    pil/.style={
           ->,
           shorten <=4pt,
           shorten >=4pt,},
    pildotted/.style={
           ->,
           shorten <=4pt,
           shorten >=4pt,
  dotted,
  },
  external/system call={pdflatex \tikzexternalcheckshellescape 
                                        -halt-on-error
                                        -interaction=batchmode 
                                        -jobname "\image" "\texsource"
                                        && pdftops -eps "\image.pdf"}
}

\usepackage{xurl}

\endlocaldefs

\begin{document}

\begin{frontmatter}
\title{Estimation for Multistate Models Subject to Reporting Delays and Incomplete Event Adjudication With Application to Disability Insurance}
\runtitle{Estimation with Reporting Delays and Incomplete Event Adjudication}

\begin{aug}
\author[A]{\fnms{Kristian}~\snm{Buchardt}\ead[label=e1]{kristian@buchardt.net}},
\author[B]{\fnms{Christian}~\snm{Furrer}\ead[label=e2]{furrer@math.ku.dk}}
\and
\author[B]{\fnms{Oliver}~\snm{Sandqvist}\ead[label=e3]{oliver.s@math.ku.dk}}
\address[A]{AP Pension\printead[presep={,\ }]{e1}}
\address[B]{Department of Mathematical Sciences, University of Copenhagen\printead[presep={,\ }]{e2,e3}}
\end{aug}

\begin{abstract}

Accurate forecasting of an insurer's outstanding liabilities is vital for the solvency of insurance companies and the financial stability of the insurance sector. For health and disability insurance, the liabilities are intimately linked with the biometric event history of the insured. Complete observation of event histories is often impossible due to sampling effects such as right-censoring and left-truncation, but also due to reporting delays and incomplete event adjudication. In this paper, we develop a parametric two-step M-estimation method that takes the aforementioned effects into account, treating the latter two as partially exogenous. The approach is valid under weak assumptions and allows for complicated dependencies between the event history, reporting delays, and adjudication while remaining relatively simple to implement. The estimation approach has desirable properties which are demonstrated by theoretical results and numerical experiments.

In the application, we introduce and consider a large portfolio of disability insurance policies. We find that properly accounting for the sampling effects has a large impact on the number of disabilities and reactivations that an insurer would forecast, allowing for a more accurate assessment of the insurer's liabilities and improved risk management.
\end{abstract}


\end{frontmatter}


\section{Introduction}
\label{sec:intro}

Health and disability insurance is a fundamental pillar of modern healthcare systems, providing individuals with financial protection against the costs of medical treatment and the loss of income due to injury or illness. In 2023, the health insurance industry of the United States collected more than \$1 trillion in premiums and paid out hundreds of billions of dollars in medical claims according to the~\citet{NAIC:2023}. For insurers, managing future liabilities such as outstanding disability benefits is crucial not only for their own financial sustainability but also for the stability of the entire health insurance market. Accurate forecasting of these liabilities, so-called reserving, is therefore vital in ensuring that insurers can meet their obligations and maintain solvency.

Insurance benefits for health and disability insurance coverages are directly tied to the insured's biometric state, making it possible to focus on modeling this underlying process. That approach is typically preferable to modeling individual payments and their intertemporal dependencies directly; this is, on the other hand, usually required for individual reserving models in non-life insurance, see for example~\citet{Yang.etal:2024}. Multistate models provide a natural and parsimonious way to model the biometric state process of an individual and are therefore the focus of this paper. By modeling individual policies, the insurer may leverage granular data to predict and forecast outstanding liabilities.

We develop estimation methods for multistate models that allow actuaries and statisticians to employ all recently generated data while accommodating the contamination induced by reporting delays and incomplete event adjudications. This assists insurers in performing timely operational adjustment and risk mitigation when faced with emerging health trends such as the major increase in claims related to mental illnesses that has been observed in recent years, see for example Section 1 of~\citet{SwissRe:2022}.

Longitudinal biometric data such as those arising from insurance policies consists of records on the occurrence and timing of certain events. The analysis of such data is typically referred to as event history analysis.
Complete observation of the event history is often rendered impossible due to sampling effects, encompassing censoring and truncation as well as situations where information about individuals under study is not up to date or correct; this may for example be due to periodic sampling, reporting delays, or incomplete event adjudication. Incomplete event adjudication refers to the situation where it is undetermined, at the time of analysis, whether some reported events satisfy predetermined criteria for being true events. Ignoring these mechanisms leads to misleading and biased analyses. Censoring and truncation are most fundamental to the field of event history analysis and have so far received the most attention. 

Our main focus is on handling reporting delays and incomplete event adjudication in the estimation of the conditional hazards of a multistate model, but we also allow for right-censoring and left-truncation. The goal is to use the available information to obtain as reliable predictions as possible. The available information consists of event dates and their corresponding event types, reporting delays, and adjudication information, but only for those events that are reported before the time of analysis.

\subsection{Relevant literature} \label{subsec:literature}

A naive way to handle reporting delays is back-censoring, see~\citet{Casper:Cook:2012}. Here, the right-censoring date is set back by some fixed amount, corresponding to only using data that is older than a given date. The approach does not introduce systematic bias if the right-censoring time is set back by an amount larger than the longest delay but is inefficient in that it may discard many observations. Since observations for later times are discarded, tail-estimates and trends may be especially impaired. This impedes forecasting and hence the ability to predict outstanding liabilities. 

More efficient approaches have been proposed for survival models, see~\citet{Hu:Tsiatis:1996} and~\citet{Van:Hubbard:1998}, and recurrent event models, see~\citet{Kalbfleisch:Lawless:1989},~\citet{Pagano.etal:1994},~\citet{Becker:Cui:1997},~\citet{Casper:Cook:2012}, and~\citet{Verbelen.etal:2022}. Most of these focus on non-parametric estimation with inverse probability of censoring weighting being the most popular approach, but other approaches include an EM-algorithm for a multinomial model in~\citet{Pagano.etal:1994}, maximum likelihood via thinning of a Poisson process in~\citet{Kalbfleisch:Lawless:1989}, and a mix of these in~\citet{Verbelen.etal:2022}. In most of the aforementioned works, the marginal mean of the counting process is modeled, except in~\citet{Kalbfleisch:Lawless:1989} and~\citet{Verbelen.etal:2022} where hazards are modeled. Finally, some works take estimation of the reporting delay distribution as the main focus and either forgo estimation of the event hazards, see \citet{Lagakos.etal:1988} and~\citet{Kalbfleisch:Lawless:1991}, or use simple methods to adjust the event rates after estimating the reporting delay distribution, see~\citet{Esbjerg.etal:1999}. 

The incorporation of reporting delays has, to the best of our knowledge, not hitherto been studied for non-competing-risks multistate models such as the illness-death models that are relevant for modeling disability insurance events. Additionally, we are the first to explore estimation under a monotone reporting assumption, which requires events to be reported in the same order in which they occurred. This assumption is natural when modeling individuals. Existing methods for estimating hazard rates in the presence of reporting delays assume that these delays are independent of the event process; this assumption only seems plausible for aggregate models whenever multiple jumps are possible. One of the primary challenges addressed by this paper is thus quantifying how reporting delays impact the hazard, assuming monotone reporting. 

For incomplete event adjudication, a simple but inefficient solution is to use only confirmed events or all unrefuted events while back-censoring. For the method to be unbiased, one has to back-censor by more than the maximal sum of the reporting and processing delays, in which case the unrefuted and confirmed events are the same. Other approaches are explored in~\citet{Cook:2000},~\citet{Cook:Kosorok:2004}, and~\citet{Bladt:Furrer:2024}, but only for single events without reporting delays. These methods also do not easily generalize to hazard estimation for multistate models since they estimate the marginal mean/distribution function of the counting process, which does not fully determine the distribution of the multistate model when multiple jumps are possible. Furthermore, in~\citet{Bladt:Furrer:2024}, the adjudication outcome distribution is exogeneously given, while in~\citet{Cook:2000} and~\citet{Cook:Kosorok:2004} it is assumed, perhaps slightly implicitly, that adjudication outcomes are conditionally independent given the covariate and event history and that the probability of confirming an event that occurred at time $s$ does not depend on information obtained after time $s$.

In addition, limited attention has been given to the estimation of the adjudication probabilities:~\citet{Cook:2000} and~\citet{Cook:Kosorok:2004} suggest logistic regression based on completed adjudications, while~\citet{Bladt:Furrer:2024} relies on expert judgments. Using only complete case data to estimate the adjudication outcome distribution can lead to biased estimates if events with shorter adjudications are systematically more or less likely to be confirmed compared to events with long adjudications. This is not an issue for the application considered by~\citet{Cook:Kosorok:2004}, where adjudications are handled by a committee that meets infrequently and reviews all unadjudicated cases at once. In that setting, the adjudication delay merely reflects the time until the next committee meeting. By contrast, in applications related to for instance disability insurance, claims are processed and adjudicated randomly and more or less continuously. Here, the adjudication time may be correlated with the adjudication outcome and decisions are not necessarily final, as insured may reapply within a statute of limitations.

Both reporting delays and incomplete event adjudication are prevalent in health and disability insurance. Disabilities only become known to the insurer when they are reported by the insured, at which point the event date also becomes known, as it determines the date from which the insured is eligible for disability benefits. Adjudication is performed primarily by the insurer, who considers whether the claim satisfies the criteria for disability benefits and determines at what stage the insured has recovered sufficiently to terminate benefits. This is often complicated by the need to obtain clinical assessments of the claimed disability. 

\subsection{Model considerations}

In order to produce accurate predictions, actuarial forecasts should be dynamic in the sense that they are updated as new information arrives. Hazard rates provide a convenient way to parameterize all relevant conditional distributions, so that insurance reserves and other conditional estimands may be calculated using well-known integral or differential equations, see~\citet{Hoem:Aalen:1978},~\citet{Moller:1993},~\citet{Norberg:2005},~\citet{Buchardt.etal:2015}, and~\citet{Adekambi:Christiansen:2017}. Consequently, this paper focuses on hazard estimation.

Baseline covariates and the event history may have a large effect on the distribution of future jumps, so we take a regression approach to incorporate these effects. Furthermore, we find parametric models advantageous for several reasons. Firstly, we need to extrapolate the estimates to times and durations that exceed the observation window. This holds both for the parameters of interest, being those related to the biometric events, and the nuisance parameters related to the reporting delay and adjudication processes. This makes the otherwise popular Kaplan--Meier, Aalen--Johansen, and Cox-type estimators less attractive. 

Secondly, feature engineering can be used in combination with existing expert knowledge to guide the parametric specification since the focus is on prediction rather than inference. This may be helpful in counteracting some of the challenging characteristics of disability insurance data, such as the relative rareness of events, the low signal-to-noise rate, and the moderate dimension of covariates. Many of these characteristics are shared with nowcasting and forecasting of epidemics, see for example~\citet{Noufaily.etal:2016},~\citet{Psotka.etal:2020}, and~\citet{Stoner.etal:2023}, which also constitutes a relevant area of application for our methods. As discussed in~\citet{Cook:2000},~\citet{Cook:Kosorok:2004}, and~\citet{Casper:Cook:2012}, another possible area of application is interim stages of clinical trials, but here semi- and non-parametric models may be preferred, given the focus is often on inference for non-parametrically identifiable marginal estimands. 

\subsection{Our contribution}

The contributions of this paper are twofold. Firstly, we develop a parametric approach that encompasses both reporting delays and incomplete event adjudication and which applies to hazard estimation for multistate models in general. The approach has many properties that are desirable for insurance applications and relaxes several assumptions from the literature. We also introduce an approximation that simplifies the estimation procedure considerably and which performs well in many situations encountered in practice. In addition, we show that our estimators are consistent and asymptotically normal under suitable regularity conditions and that they may be bootstrapped.

Secondly, the proposed methods are applied to a new data set based on a large Danish disability insurance portfolio. This application is noteworthy for two reasons. First, the data is unique in that it seems to be the first insurance data set to include information on event adjudication and the first disability or health insurance data set to contain information on reporting delays. Having access to such data and accounting for the time of analysis is essential to ensure unbiased analyses. Second, the analyses show that properly accommodating for reporting delays and incomplete event adjudication has a substantial effect on the estimated level and calendar-time dependence of the disability and reactivation hazards. This underscores the importance of our concepts and results. The impact of our contribution is further accentuated by taking the size and societal importance of the health and disability insurance sector into account.

The paper is organized as follows. Section~\ref{sec:Motivation} briefly discusses the motivating application towards disability insurance. Section~\ref{sec:Model} introduces the model. Section~\ref{sec:Estimation} discusses estimation of the parameters. Section~\ref{sec:Asymptotics} presents theoretical large-sample results. Section~\ref{sec:Numerical} contains numerical experiments, demonstrating desirable finite sample performance and stability under misspecification. The data application related to disability insurance is presented in Section~\ref{sec:Application}.

\section{Motivating Insurance Application} \label{sec:Motivation}

Multistate models are widely used in actuarial science, with transitions between states representing changes to the state of the insured. In the data application of Section~\ref{sec:Application}, the states are active (1), disabled (2), reactivated (3), and dead (4), and the transitions represent the times where the corresponding events occur. For example, if an insured construction worker falls and sustains a back injury that leaves them unable to work full time, the time of the fall would be the time of the transition from state 1 to state 2. If they recover sufficiently for the insurer to terminate disability benefits, the time where benefits were terminated would be represented by a transition from state 2 to state 3. Multistate models can generally be described by an initial state and a multivariate counting process $\mathbf{N}^\ast$ describing the transitions between the different states. State-of-the-art multistate methodology in actuarial science includes the so-called smooth semi-Markov approach, where in particular $\mathbf{N}^\ast$ is required to admit an intensity process that depends only on time, the current state, and the duration spent in the current state. Empirical evidence shows that the probability of reactivation depends greatly on the time spent as disabled, so for the disability insurance application of Section~\ref{sec:Application}, we employ such a smooth semi-Markov model rather than, for example, a Markov chain model. 

The state of the insured is crucial for reserving since its current and past values determine the benefits that the insured is eligible for, meaning that the accumulated payment process $B$ is adapted to the information generated by $\mathbf{N}^\ast$. The goal of reserving is then the prediction of future expected liabilities, which includes the estimation of the conditional mean of $B(s) - B(t)$ for $s>t$, given the information available at time $t$. In the smooth semi-Markov approach, the intensity of $B$ also only depends on time, the current state, and the duration spent in the current state. In this case, these conditional means may be computed by solving the system of differential equations given in Theorem 3.1 of~\citet{Buchardt.etal:2015} or Corollary 7.2 of~\citet{Adekambi:Christiansen:2017} once the hazards have been estimated. For more details on multistate modeling in actuarial science, we refer to Sections~1.1 and~1.2 in~\citet{Furrer:2020}.

In the data application of Section~\ref{sec:Application}, we introduce and study a new data set related to disability insurance and estimate the relevant hazard rates. Considering a short coverage period such that multiple distinct disabilities for an individual are unlikely, it is common to adopt a four-state semi-Markov model corresponding to an illness-death model with a separate reactivated state. Let $I^\ast_j(t)$ be the indicator that the state of the insured is $j$ at time $t$ and let $D^\ast(t)$ be the duration spent in the current state at time $t$. An important example of a disability insurance payment process is then
\begin{align*}
B(\mathrm{d}t) = 1\{t-D^\ast(t) \leq q\} \{ w \leq D^\ast(t)\}  I^\ast_2(t) b \, \mathrm{d}t,  
\end{align*}
which corresponds to a benefit rate $b$ during disability,  but subject to a coverage period of length $q$ and a waiting period of duration $w$. Since contract features such as the benefit rate, the coverage period, and the waiting period may differ between insured, and due to left-truncation and right-censoring effects, hazard estimation is recommendable. Left-truncation corresponds to random entries into the insurance portfolio from some reference population, and right-censoring  corresponds to potentially random exits, which is commonly state-dependent as insured more regularly leave when active or reactivated. However, as explained in the introduction, transitions from state $1$ to $2$ (disabilities) are affected by reporting delays, and transitions from state $1$ to state $2$ (disabilities) and from $2$ to $3$ (reactivations) are subject to adjudication. If this is not appropriately addressed, the insurer risks systematically over- or underestimating its liabilities, which leads to worsened risk management.

The four-state semi-Markov model considered in Section~\ref{sec:Application} is perhaps the simplest actuarial mulitstate model imaginable for which both delays and adjudications are relevant. In practice, even just for disability insurance, more complicated models are required to capture multiple distinct disabilities, different types of disabilities, and the interplay with other coverages such as health insurance. Consequently, there is a need for hazard estimation methodology that can handle not only left-truncation and right-censoring, but also delays and adjudications, for mulitstate models in general.

\section{Model Specification} \label{sec:Model}

We consider a time-continuous multistate model with states $1, \ldots, J$, which may be represented by a marked point process $(T^\ast_m,Y^\ast_m)_{m\geq1}$, where $T^\ast_m \in (0,\infty]$ are jump-times and $Y^\ast_m \in \{1,\ldots,J\}$ are jump-marks.  Let $\mathbf{J}^\ast(t)$ denote the time $t$ event history, containing the events that have occurred up until and including time $t$. The multistate model may equivalently be represented by a multivariate counting process $\mathbf{N}^\ast$ with components $N^\ast_{jk}$, so that $N^\ast_{jk}(t)$ counts the number of transitions from state $j$ to state $k$ in the interval $[0,t]$. Let $\mathbf{X}$ be baseline covariates, including for instance the initial state and age of the subject, taking values in a subset of $\mathbb{R}^p$. We suppose that $N^\ast_{jk}$ satisfies a general intensity model
\begin{align}\label{eq:mu_ast}
\mathbb{E}_{\boldsymbol \theta}  \{ N^\ast_{jk}(\mathrm{d}t) \mid \mathcal{F}^\ast_{t-} \} = \lambda^\ast_{jk}(t; \mathbf{X}, \boldsymbol\theta) \, \mathrm{d}t = I^\ast_j(t)\mu^\ast_{jk}\{ t; \mathbf{J}^\ast(t-), \mathbf{X}, \boldsymbol\theta \} \, \mathrm{d}t,
\end{align}
where $\mathbb{E}_{\boldsymbol \theta}$ denotes expectation under the distribution with parameter $\boldsymbol\theta \in \Theta \subseteq \mathbb{R}^d$, $\mathcal{F}^\ast$ is the filtration generated by $\mathbf{N}^\ast$ and $\mathbf{X}$, and $I^\ast_j(t)$ is the indicator that the the last mark in $\mathbf{J}^\ast(t-)$ is $j$ or the initial state is $j$ if $\mathbf{J}^\ast(t-)$ is empty.  If the intensity only depends on $\mathbf{J}^\ast(t-)$ through the latest jump-time and jump-mark, then one obtains a smooth semi-Markov model as described in Section~\ref{sec:Motivation}. Assume that the true parameter $\boldsymbol\theta_0$ belongs to the parameter set $\Theta$. We seek to estimate $\boldsymbol\theta_0$, but subject to various forms of missingness and contamination. First, the observation of jumps is subject to delays. Second, not all jump-times, jump-marks, and delays are observed due to random left-truncation and right-censoring. Finally, some reported jumps are annulled due to event adjudication. In the following, we describe these mechanisms in more detail.

Let $\eta>0$ be the deterministic time where the statistical analysis is conducted. Denote by $U^\ast_m \in [0,\infty)$ the random reporting delay associated with the $m$'th jump. We form a new thinned marked point process $(T_m, Y_m)_{m \geq 1}$ from $(T_m^\ast,Y_m^\ast)_{m \geq 1}$ by deleting all jumps with $T^\ast_m + U^\ast_m > \eta$. Let $(U_m)_{m \geq 1}$, $\mathbf{J}$, and $N_{jk}$ be respectively the reporting delays, event history and counting processes associated with $(T_m, Y_m)_{m \geq 1}$. Throughout, we make the following assumption:
\begin{assumption} \label{AssumptionMonotoneReporting}
    Reporting is monotone, meaning $T^\ast_{m} + U^\ast_{m} \leq T^\ast_{m+1} + U^\ast_{m+1}$.
\end{assumption}
\noindent Assumption~\ref{AssumptionMonotoneReporting} is closely related to the concept of monotone missingness from the missing data literature, confer with~\citet{Little:2021}. Dropping this assumption would lead to a substantial increase in complexity, since even if $T_m < \infty$, there might still be a disagreement between $\mathbf{J}^\ast(T_m)$ and $\mathbf{J}(T_m)$.  This assumption is natural when modeling individual subjects. Thinning has hitherto only been used for aggregate models when this assumption is replaced by an independence assumption. The true distribution function of $U_m^\ast$ given $\mathbf{J}^\ast(T^\ast_m)$ and $\mathbf{X}$ is denoted $t \mapsto \textnormal{pr}_{U}\{ t ; \mathbf{J}^\ast(T^\ast_m),\mathbf{X}, \mathbf{f}_0 \}$ for a parameter $\mathbf{f}_0$. That the components of $(U^\ast_m)_{m \geq 1}$ may have different distributions is captured by $\mathbf{J}^\ast(T^\ast_m)$ in the conditioning. 

Observation of $\mathbf{N}$ is, as mentioned, also subject to left-truncation and right-censoring. To be precise, there is a random right-censoring time $C \leq \eta$, which could for example be the drop-out time of the subject, and a random left-truncation time $V < C$ where the subject enters the study if some event $A$ has occurred prior to $V$. As an example, $V$ could be when the subject enters the insurance portfolio, and the event $A$ could be that the subject is still alive at that time. A comparable situation would be one where subjects not satisfying $A$ or with $V > \eta$ were also considered part of the study but with no exposure. For example, if one had access to panel data, the current formulation would filter away subjects that are not in the portfolio during the observation window while the alternative formulation would keep the observations but let them have no exposure. In our applications, the event $A$ is redundant as it is contained in the relevant filtration, but we allow for this extra generality which may be useful in other applications, confer with Section III.3 of~\citet{Andersen.etal:1993}. 

Let $(T^c_m,Y^c_m)_{m \geq 1}$ be $(T_m,Y_m)_{m \geq 1}$ but where jumps outside of $(V,C]$ are deleted and let $(U^c_m)_{m \geq 1}$, $N^c_{jk}$, and $\mathbf{J}^c$ be the corresponding reporting delays, counting processes, and event history, respectively. We introduce the filtrations $\mathcal{F}_t = \sigma \{ \mathbf{X}, (T_m, Y_m) : T_m \leq t  \}$  and  $\mathcal{F}^c_t = \sigma \{\mathbf{X}, \mathcal{G}, C \wedge t, (T^c_m, Y^c_m) : T^c_m \leq t \}$, where $\mathcal{G}$ is left-truncation information satisfying $A \in \mathcal{G}$ and $V$ is $\mathcal{G}$-measurable; the delays are not included in these filtrations as the influence of prior delays on the intensities is not of interest when estimation of $\mu^\ast_{jk}$ is the final goal. Also, these filtrations are not directly observable (for instance, $T_m$ is reported at time $T_m+U_m$, but already enters in the filtrations at time $T_m$), but they play a technical role in allowing us to formulate a criterion function from which to construct estimators. In our application, we take $\mathcal{G}=\sigma(V)$. The more general notion of $\mathcal{G}$ is included as it might be useful in other situations. e.g.\ if subjects could enter the portfolio as disabled, in which case we would need $\mathcal{G}$ to contain the duration as disabled in order to be able to compute the hazards.

Note that we allow for jumps to be reported between the right-censoring time $C$ and the time of analysis $\eta$; this choice is appropriate for actuarial applications where subjects may leave the portfolio and afterwards report claims for events that occurred while they were still in the portfolio. In other applications, such as for medical trials, thinning according to $T^\ast_m + U^\ast_m > C$ might be more appropriate, as patients may be completely lost to follow-up at $C$. We stress that the results of this paper then still apply upon changing $\eta$ to $C$ in the reporting delay distribution, but under the additional assumption that $C$ is independent of $(T^\ast_m,Y^\ast_m)_{m\geq 1}$. Furthermore, in some applications it might be relevant to also allow for periodic, but immediate, ascertainment of the current state. This has been explored for a survival model in~\citet{Hu:Tsiatis:1996}. We do not pursue this extension. 

The following assumption, which we make throughout, is comparable with the assumption of independent filtering, confer with Section~III.4 in~\citet{Andersen.etal:1993}.
\begin{assumption} \label{assumption:IndepCensoring}
    Left-truncation and right-censoring are independent in relation to the marked point process $(T_m,Y_m)_{m \geq 1}$ in the sense that
\begin{align*}
\mathbb{E}_{\boldsymbol\theta,\mathbf{f}} \{ N^c_{jk}(\mathrm{d}t) \mid  \mathcal{F}^c_{t-};A \} = 1_{(V,C]}(t)\mathbb{E}_{\boldsymbol\theta,\mathbf{f}} \{ N_{jk}(\mathrm{d}t) \mid  \mathcal{F}_{t-} \}
\end{align*}
for $(\boldsymbol\theta,\mathbf{f}) \in \Theta \times \mathbb{F}$, where the notation on the left-hand side signifies conditioning on $\mathcal{F}^c_{t-}$ and $A$. Similarly, left-truncation and right-censoring are independent in relation to the reporting delays in the sense that $\mathbb{P}_{\mathbf{f}}(U^c_m \leq t \mid \mathcal{F}^c_{T_m^c};A) = \mathbb{P}_{\mathbf{f}}(U_\ell \leq t \mid \mathcal{F}_{T_\ell})$ where $\ell$ is the index satisfying $T_\ell=T^c_m$. Here $\mathbb{P}_{\mathbf{f}}$ denotes the probability under the distribution with parameter $\mathbf{f}$.
\end{assumption}

The third and final mechanism relates to incomplete event adjudication. The idea is to think of $(T^c_m,Y^c_m)_{m \geq 1}$ as obtained by thinning another marked point process $(\tilde{T}_m,\tilde{Y}_m)_{m \geq 1}$ accompanied by delays $(\tilde{U}_m)_{m\geq1}$ satisfying $\tilde{T}_m + \tilde{U}_m \leq \eta$. To be precise, if $\xi_m \in \{0,1\}$ is the adjudication outcome of the $m$'th event, then $(T^c_m,Y^c_m)_{m \geq 1}$ is formed from $(\tilde{T}_m,\tilde{Y}_m)_{m \geq 1}$ by deleting all jumps with $\xi_m = 0$. Incomplete event adjudication thus becomes a missing data problem for the adjudication outcomes $(\xi_m)_{m \geq 1}$. The full available information is thus
\begin{align*}
    \mathcal{H}^{\textnormal{obs}}_t &= \sigma \{\mathbf{X},\mathcal{G},C \wedge t,(\Tilde{T}_m,\Tilde{Y}_m,\Tilde{U}_m, \mathcal{A}_{m,t}) : \Tilde{T}_m+\Tilde{U}_m \leq t \},
\end{align*}
where ${\mathcal{A}}_{m,t}$ is a filtration representing adjudication information of the $m$'th reported jump, such that when $\sigma_m$ is the time where $\xi_m$ becomes known, then $\sigma_m$ is an $\mathcal{A}_{m,t}$-stopping time and $\xi_m$ is $\mathcal{A}_{m,\sigma_m}$-measurable. Apart from $\sigma_m$ and $\xi_m$, the filtration $\mathcal{A}_{m,t}$ may contain additional covariates that affect adjudication probabilities.  In Section~\ref{sec:est_h}, we propose to let  ${\mathcal{A}}_{m,t}$ be the information generated by a multistate model representing the adjudication process and to let $\xi_m$ be determined by which state the multistate model becomes absorbed in, such that $\sigma_m$ is the absorption time.

The following assumption serves to reduce the modeling task:
\begin{assumption} \label{assumption:oneAdj}
    Only the most recent jump of $(\Tilde{T}_m,\Tilde{Y}_m)_{m \geq 1}$ may be unadjudicated at any given time, meaning $\sigma_m \leq \Tilde{T}_{m+1}+\Tilde{U}_{m+1}$.
\end{assumption}
\noindent No mathematical issues would arise from removing this assumption, but one would then need to find suitable estimators for the joint distribution of the adjudication outcomes instead of for single outcomes. 

We also introduce a restricted (non-monotone) information $\mathcal{H}_t \subset \mathcal{H}^{\textnormal{obs}}_t$, which we use as the conditioning information in the distribution of the adjudication outcomes. The point is that it may sometimes be convenient to let the adjudication hazards depend on less information than all available information. It is important that the restricted information contains the confirmed jumps so that they are not treated as missing values, but which additional information one lets the adjudication probabilities depend on is optional. As will be seen in Section~\ref{sec:Estimation}, $\mathcal{H}_t$ is needed for all $t \leq \eta$ to estimate the adjudication model, while $\mathcal{H}_\eta$ is sufficient to estimate the reporting delays and biometric hazards. Denote by $\langle t \rangle$ the number reported jumps of $(\Tilde{T}_m,\Tilde{Y}_m)_{m \geq 1}$ before time $C \wedge t$. We find a natural choice of $\mathcal{H}_t$ to be
\begin{align*}
    \mathcal{H}_t &= \sigma \{\mathbf{X}, \mathcal{G}, C \wedge t ,(\Tilde{T}_m,\Tilde{Y}_m),(\Tilde{T}_{\langle t \rangle},\Tilde{Y}_{\langle t \rangle},\Tilde{U}_{\langle t \rangle}, \mathcal{A}_{\langle t \rangle,t}) : \Tilde{T}_m+\Tilde{U}_m \leq t, \xi_m = 1, m < \langle t \rangle \}
\end{align*} 
consisting of the confirmed jumps by time $t$ as well as the most recent reported jump with its reporting delay and adjudication information; information about prior reporting delays and adjudication processes is removed. This is the choice of $\mathcal{H}$ that is used in the application, allowing us to avoid modeling how the reactivation adjudication hazards would be affected by adjudication information related to the disability event, for example how long it took for the disability to be confirmed. We suppose that the distribution of the adjudication outcome $\xi_{\langle t\rangle}$ given $\mathcal{H}_t$ is characterized by a parameter $\mathbf{g}_0$. The parameter set for $(\mathbf{f},\mathbf{g})$ is assumed to be a product space $\mathbb{F} \times \mathbb{G}$ containing $(\mathbf{f}_0,\mathbf{g}_0)$. 

It is worth noting that we place no distributional assumptions on  $(\tilde{T}_m,\tilde{Y}_m)_{m \geq 1}$ apart from those inherited from $(T^c_m,Y^c_m)_{m \geq 1}$. Likewise, we do not want to impose criteria for how our process of interest $(T^\ast_m,Y^\ast_m)_{m \geq 1}$ would be affected by subsampling on adjudication information or reporting delays. We therefore choose to treat incomplete event adjudication and reporting delays as partially exogenous mechanisms. This naturally leads to a two-step method, where the first step concerns incomplete event adjudication and reporting delays, while the second step involves estimation of $\boldsymbol\theta_0$. The method is outlined in Section~\ref{sec:Estimation}, and the main theoretical results, including weak consistency and asymptotic normality, are given in Section~\ref{sec:Asymptotics}. An alternative approach described in Section C of the Supplementary material (\citet{Buchardt.etal:2025}) incorporates the delays in the second step; this improves efficiency, but comes at the cost of additional independence assumptions. 

\section{Estimation} \label{sec:Estimation}

\subsection{General considerations} \label{subsec:General}

For our two-step method, we propose a two-step M-estimation procedure. Results on the asymptotics of such procedures may be found in~\citet{Murphy:Topel:1985} and~\citet{Newey:McFadden:1994} for the parametric case and~\citet{Ichimura:Lee:2010},~\citet{Kristensen:Salanie:2017}, and~\citet{Delsol:Van:2020} for the semiparametric case. We treat $\mathbf{f}$ and $\mathbf{g}$ as nuisance parameters and $\boldsymbol\theta$ as the parameter of interest. Discussion of concrete estimators is postponed to Section~\ref{sec:est_h}.

For notational simplicity, we henceforward suppress the baseline covariates and use the symbol $\bigcdot$ to signify summation over the corresponding index. Disregarding for a moment incomplete event adjudication, in the sense that we work under a filtration where the adjudication outcomes for all reported jumps are known, and fixing the parameter of the distribution of the reporting delays at $\mathbf{f}$, we obtain from Assumption~\ref{assumption:IndepCensoring} and Section~II.7.3 in~\citet{Andersen.etal:1993} the partial likelihood $L(\boldsymbol\theta;\mathbf{f})$ for one subject under $\mathcal{F}^c_\eta$ as the following product integral
\begin{align*}
L(\boldsymbol\theta;\mathbf{f}) = \Prodi_{t = V}^{C}
\prod_{j = 1}^J \bigg\{ \prod_{k \neq j} {\Lambda_{jk}(\mathrm{d}t; \boldsymbol\theta, \mathbf{f})}^{\Delta N_{jk}(t)} \bigg\} \bigg\{ (1 - {\Lambda_{j\bigcdot}(\mathrm{d}t; \boldsymbol\theta, \mathbf{f}))}^{1 - \Delta N_{j\bigcdot}(t)}\bigg \},
\end{align*}
where $\Lambda_{jk}$ is the compensator of $N_{jk}$ i.e.\ $\Lambda_{jk}(\mathrm{d}t; \boldsymbol\theta, \mathbf{f}) = \mathbb{E}_{\boldsymbol\theta, \mathbf{f}} \{N_{jk}(\mathrm{d}t) \mid  \mathcal{F}_{t-}\}$.
 Denote the corresponding log-likelihood by $\ell(\boldsymbol\theta;\mathbf{f}) = \log L(\boldsymbol\theta;\mathbf{f})$. We henceforth assume that $\mathbb{E}[\ell(\boldsymbol\theta;\mathbf{f})]$ is uniquely maximized in $(\boldsymbol\theta_0,\mathbf{f}_0)$ which enables the use of M-estimation techniques. This holds under weak assumptions. For example, Section~II.7.2 in~\citet{Andersen.etal:1993} implies, under some integrability and smoothness conditions, that $(\boldsymbol\theta_0,\mathbf{f}_0)$ is a unique maximum of $\mathbb{E}[\ell(\boldsymbol\theta;\mathbf{f})]$ as long as the compensators are non-constant in each coordinate of $(\boldsymbol\theta,\mathbf{f})$. 

Seeking a tractable expression of $L$ in terms of~\eqref{eq:mu_ast}, we introduce the weighted hazard
\begin{align*}
\nu_{jk}\{t; \mathbf{J}(t-), \boldsymbol\theta, \mathbf{f}\}
=
\mu^\ast_{jk}\{t; \mathbf{J}(t-), \boldsymbol\theta \} \textnormal{pr}_{U}\big[\eta-t ; \{\mathbf{J}(t-),t,k\}, \mathbf{f} \big]
\end{align*}
with $\{\mathbf{J}(t-),t,k\}$ being the event history containing the jumps $\mathbf{J}(t-)$ and a jump to state $k$ at time $t$. Let $D(t)$ denote the duration since the most recent jump of $(T_{m},Y_m)_{m \geq 1}$ at time $t$, let $I_j(t)$ be the indicator that the last mark in $\mathbf{J}(t-)$ is $j$ or the initial state is $j$ if $\mathbf{J}(t-)$ is empty. Introduce the survival probability in the current state $P^\ast \{t; \mathbf{J}(t-), \boldsymbol\theta \}$ satisfying
\begin{align*}
I_j(t)P^\ast \{t; \mathbf{J}(t-), \boldsymbol\theta \} = I_j(t)\exp\!\bigg[-\int_{t-D(t-) }^t \mu^\ast_{j\bigcdot}\{s; \mathbf{J}(s-), \boldsymbol\theta \} \, \mathrm{d}s \bigg].
\end{align*}
In Section B of the Supplementary material, we derive the following result.
\begin{lemma}\label{lemma:hazard_link}
It holds that $\Lambda_{jk}(\mathrm{d}t; \boldsymbol\theta, \mathbf{f}) = I_j(t)\mu_{jk} \{t; \mathbf{J}(t-), \boldsymbol\theta, \mathbf{f} \} \, \mathrm{d}t$ with
\begin{align*}
\mu_{jk}\{t; \mathbf{J}(t-), \boldsymbol\theta, \mathbf{f}\}
=
\gamma_j(t) \times 
\nu_{jk}\{t; \mathbf{J}(t-), \boldsymbol\theta, \mathbf{f}\}
\end{align*}
where $\gamma_j(t)$ equals
\begin{align*}
    \frac{1-\int_{t-D(t-)}^t P^\ast \{s; \mathbf{J}(s-), \boldsymbol\theta \} \mu^\ast_{j\bigcdot} \{s; \mathbf{J}(s-), \boldsymbol\theta\} \, \mathrm{d}s}
{\textnormal{pr}_U[\eta - \{t-D(t-)\}; \mathbf{J}\{t-D(t-)\}, \mathbf{f}]-\int_{t-D(t-)}^t P^\ast\{s; \mathbf{J}(s-), \boldsymbol\theta\} \nu_{j\bigcdot}\{s; \mathbf{J}(s-), \boldsymbol\theta, \mathbf{f}\} \, \mathrm{d}s}.
\end{align*}
\end{lemma}

\noindent We now return to the problem that $\ell(\boldsymbol\theta;\mathbf{f})$ is actually not computable using the available information $\mathcal{H}^{\textnormal{obs}}_\eta$ due to incomplete event adjudication. Furthermore, estimation of $\mathbf{f}$ based on $\ell(\boldsymbol\theta;\mathbf{f})$ is inefficient since it only utilizes $(T_m,Y_m)_{m \geq 1}$ and not $(U_m)_{m \geq 1}$. We therefore propose two-step M-estimation. Denote by $\mathbf{Z}_{i}^{\textnormal{obs}}$ and $\mathbf{Z}_{i}$ $(i=1,\ldots,n)$ the subject $i$ outcomes of $\mathcal{H}^\textnormal{obs}_\eta$ and $\mathcal{H}_\eta$ respectively, and let $\mathbf{Z}^{\textnormal{obs}}$ and $\mathbf{Z}$ be corresponding generic outcomes. Assume $\mathbf{Z}_i^{\textnormal{obs}}$ are independent and identically distributed, which implies that the same holds for $\mathbf{Z}_i$. We specify the objective function to be the imputed likelihood $\ell(\mathbf{Z},\boldsymbol\theta;\mathbf{f},\mathbf{g}) =  \mathbb{E}_{\mathbf{g}}\{\ell(\boldsymbol\theta;\mathbf{f}) \mid \mathbf{Z}\}$. For the relation to imputation, see Section A of the Supplementary material. Our two-step procedure is:
\begin{enumerate}
\item Estimate $(\mathbf{f}_0, \mathbf{g}_0)$ by $(\hat{\mathbf{f}}_n, \hat{\mathbf{g}}_n)$ using suitable estimators based on $(\mathbf{Z}^\textnormal{obs}_i)_{i=1}^n$.
\item Estimate $\boldsymbol\theta_0$ by $\hat{\boldsymbol\theta}_n=
\argmax_{\boldsymbol\theta} \sum_{i=1}^n \ell(\mathbf{Z}_{i},\boldsymbol\theta; \hat{\mathbf{f}}_n,\hat{\mathbf{g}}_n)$.
\end{enumerate}
This is motivated by the observation $\mathbb{E}\{\ell(\mathbf{Z},\boldsymbol\theta;\mathbf{f}_0,\mathbf{g}_0)\} = \mathbb{E}\{\ell(\boldsymbol\theta;\mathbf{f}_0)\}$, which is uniquely maximized in $\boldsymbol\theta = \boldsymbol\theta_0$. Recall that $\mathcal{H}^{\textnormal{obs}}_\eta$ and $\mathcal{H}_\eta$ are such that only $\xi_{\langle \eta \rangle}$ can be missing in $\ell(\boldsymbol\theta;\mathbf{f})$, and that we wish to use adjudication probabilities based on the information $\mathcal{H}_\eta$ instead of $\mathcal{H}^{\textnormal{obs}}_\eta$, leading us to condition on $\mathbf{Z}$ instead of $\mathbf{Z}^\textnormal{obs}$ in the objective function. We therefore let $w(1,\mathbf{Z};\mathbf{g}) = \mathbb{E}_\mathbf{g}\{\xi_{\langle \eta \rangle} \mid \mathbf{Z}\}$ and $w(0,\mathbf{Z};\mathbf{g})=1-w(1,\mathbf{Z};\mathbf{g})$. Write $\ell(\boldsymbol\theta;\mathbf{f}) = \ell(\boldsymbol\theta;\mathbf{f},\xi_{\langle \eta \rangle})$, where the last argument signifies whether the jump $(\Tilde{T}_{\langle \eta \rangle},\Tilde{Y}_{\langle \eta \rangle})$ is included in the likelihood or not. We thus have 
\begin{align} \label{eq:FinalLikelihood}
    \ell(\mathbf{Z},\boldsymbol\theta;\mathbf{f},\mathbf{g})=w( 0,\mathbf{Z}; \mathbf{g}) \ell(\boldsymbol\theta; \mathbf{f}, 0)+w( 1,\mathbf{Z}; \mathbf{g}) \ell(\boldsymbol\theta; \mathbf{f}, 1).
\end{align}
In Section~\ref{sec:est_h}, we outline the first step, while Section~\ref{sec:est_theta} is dedicated to the second. The derivation of asymptotic properties is postponed to Section~\ref{sec:Asymptotics}.

\subsection{Estimation of $\mathbf{f}$ and $\mathbf{g}$}\label{sec:est_h}

Observe that Assumption~\ref{AssumptionMonotoneReporting} implies $\mathbb{P}_{\mathbf{f}}\{U_m \leq t \mid \mathbf{J}(T_m)\} = \textnormal{pr}_U\{t;\mathbf{J}(T_m),\mathbf{f}\}/\textnormal{pr}_U\{\eta - T_m;\mathbf{J}(T_m),\mathbf{f}\}$ 
on the event $(T_m < \infty)$ for $t \leq \eta-T_m$. Note that one is not conditioning on $\mathbf{J}(T_m)$ on the right-hand side but rather inputting the values into the regular conditional probability $\textnormal{pr}_U$. Thus, the parameter $\mathbf{f}$ could be estimated if $(T_m,Y_m)_{m \geq 1}$ rather than $(\Tilde{T}_m,\Tilde{Y}_m)_{m \geq 1}$ were available. We therefore also employ two-step M-estimation for the nuisance parameters, first constructing $\hat{\mathbf{g}}_n$ and subsequently using an imputed likelihood to construct $\hat{\mathbf{f}}_n$. Practical aspects of the implementation are discussed in Section D of the Supplementary material. 

 We first consider estimation of $\mathbf{g}$. Since there is time dependence in adjudication processing, and since we see no need for infinite-variation processes, we choose to model the adjudication information $\mathcal{A}_{m,t}$ of the $m$'th reported jump as being generated by a marked point process $(\tau_{m,\ell},Y_{m,\ell} )_{\ell \geq 1}$ with finite state space $1,\dots,K$. The time $\sigma_m$ is specified as a first hitting time of a subset of the states $1,\dots,K$, and the value of $\xi_m$ implied by hitting one of these states is non-random. The corresponding counting processes are denoted $N_{m,jk}$ and the indicator that the process is in state $j$ at time $t-$ is denoted $I_{m,j}(t)$. This multistate model starts when the corresponding jump is reported and ends when the next jump is reported, meaning $\Tilde{T}_m+\Tilde{U}_m < \tau_{m,\ell} \leq \Tilde{T}_{m+1}+\Tilde{U}_{m+1}$ for all $\ell$. We choose a log-likelihood type objective function for $\mathbf{g}$ and parameterize via conditional hazards $\omega_{jk}(t;\mathcal{H}_{t-},\mathbf{g})$ given $\mathcal{H}_{t-}$, that is
\begin{align*}
\ell_\xi(\mathbf{Z}^\textnormal{obs},\mathbf{g}) = \sum_{m=1}^{\langle \eta \rangle} \sum_{j,k = 1}^K \int_0^\eta \log \omega_{jk}(t;\mathcal{H}_{t-},\mathbf{g}) \: N_{m,jk}(\diff t)-\int_0^\eta I_{m,j}(t)\omega_{jk}(t;\mathcal{H}_{t-},\mathbf{g}) \diff t.
\end{align*}
Note that the log-likelihood depends on $\mathcal{H}_t$ for all $t \leq \eta$; this information is contained in $\mathcal{H}^\textnormal{obs}_\eta$ but not in $\mathcal{H}_\eta$. Estimation now proceeds by maximizing $\sum_{i=1}^n \ell_\xi(\mathbf{Z}_i^\textnormal{obs},\mathbf{g})$ with respect to $\mathbf{g}$ over a suitable parametric class. This maximization is standard. Other modeling choices for $\mathbf{g}$ are possible, and while our proofs of the asymptotic properties utilize this specific choice, they can easily be adapted to a wider range of models.
 
We now consider estimation of $\mathbf{f}$. Disregarding incomplete event adjudication for a brief moment, Assumption~\ref{assumption:IndepCensoring} implies that the log-likelihood for the reporting delays becomes $$\ell_U(\mathbf{f}) = \sum_{m=1}^{N^c_{\bigcdot \bigcdot}(\eta)} \log \left[ \textnormal{pr}_{U}\{\diff U^c_m ; \mathbf{J}^c(T^c_m), \mathbf{f}\} / \textnormal{pr}_{U}\{\eta-T^c_m ; \mathbf{J}^c(T^c_m), \mathbf{f}\} \right].$$ 
Define the reverse time hazard $\alpha\{s;\mathbf{J}^c(T^c_m),\mathbf{f}\} \diff s = \textnormal{pr}_{U}\{\diff s; \mathbf{J}^c(T^c_m), \mathbf{f}\}/\textnormal{pr}_{U}\{s ; \mathbf{J}^c(T^c_m), \mathbf{f}\}$ similarly to~\citet{Kalbfleisch:Lawless:1991} and write the distribution function in terms of the reverse time hazard 
$\textnormal{pr}_U\{t;\mathbf{J}^c(T^c_m),\mathbf{f}\} = \exp\left[ -\int_t^\infty \alpha\{s;\mathbf{J}^c(T^c_m),\mathbf{f}\} \diff s \right]$. Then
\begin{align*}
\ell_U(\mathbf{f}) = \sum_{m=1}^{N^c_{\bigcdot \bigcdot}(\eta)} \log \alpha\{U^c_m;\mathbf{J}^c(T^c_m),\mathbf{f}\}  - \int_{U^c_m}^{\eta-T^c_m} \alpha\{s;\mathbf{J}^c(T^c_m),\mathbf{f}\} \diff s.
\end{align*}
Parameterizing in terms of reverse time hazards maintains the form of the likelihood under right-truncation similarly to how hazards maintain the form of the likelihood under right-censoring.
Returning to incomplete event adjudication, we analogously to $\ell(\boldsymbol\theta;\mathbf{f},\xi_{\langle \eta \rangle})$ write $\ell_U(\mathbf{f}) = \ell_U(\mathbf{f};\xi_{\langle \eta \rangle})$ and let the objective function for $\mathbf{f}$ be the imputed likelihood $\ell_U(\mathbf{Z},\mathbf{f};\mathbf{g}) = \mathbb{E}_{\mathbf{g}}\{\ell_U(\mathbf{f}) \mid \mathbf{Z}\} =  w( 0,\mathbf{Z}; \mathbf{g}) \ell_U(\mathbf{f};0)+w( 1,\mathbf{Z}; \mathbf{g}) \ell_U(\mathbf{f};1)$. For a given estimator $\hat{\mathbf{g}}_n$, we let $\hat{\mathbf{f}}_n$ be the maximizer of $\sum_{i=1}^n \ell_U(\mathbf{Z}_i,\mathbf{f};\hat{\mathbf{g}}_n)$ with respect to $\mathbf{f}$ over a suitable parametric class. 

\subsection{Estimation of $\boldsymbol\theta$}\label{sec:est_theta}

The estimator $\hat{\boldsymbol\theta}_n$ was described in Section~\ref{subsec:General}. Evaluating $\ell(\boldsymbol\theta; \mathbf{f})$ is costly due to the repeated numerical integration required owing to $\gamma_j$. This would be further compounded when evaluating the score, and we hence recommend derivative-free maximization 
for example via quasi-Newton methods. For these algorithms to perform well, one needs a good starting value for $\boldsymbol\theta$. We here describe one approach to overcoming this problem.

Define $\check{\gamma}_j$ by replacing $\textnormal{pr}_U[\eta - \{t-D(t-)\}; \mathbf{J}\{t-D(t-)\},\mathbf{f}]$ with $1$ in the expression for $\gamma_j$, effectively ignoring the reporting delay of the previous jump, and let $\check{\mu}_{jk} = \check{\gamma}_j \times \nu_{jk}$. Note that the approximation is exact if there is never a reporting delay for the prior jump or if there has been no prior jump. Note $I_j(t) P^\ast\{t; \mathbf{J}(t-), \boldsymbol\theta\} \leq \check{\gamma}_j(t) \leq 1$ so $\check{\mu}_{jk} \approx \nu_{jk}$. Furthermore, under additional smoothness, we show in Section B of the Supplementary material that $\check{\mu}_{jk}\{t;\mathbf{J}(t-),\boldsymbol\theta,\mathbf{f}\} = \nu_{jk}\{t; \mathbf{J}(t-), \boldsymbol\theta, \mathbf{f}\}+\textnormal{err}(t)$ with $\vert\textnormal{err}(t)\vert \leq D(t-) \times \sup_{s \in [t-D(t-),t]} \mu^\ast_{ j \bigcdot }\{s; \mathbf{J}(s-),\boldsymbol\theta\}^2$, so the error is a second order term. We call the approximation $\mu_{jk} \approx \nu_{jk}$ a Poisson approximation, since the approximate hazard $\nu_{jk}$ has a similar form to the one found in the Poisson process setup of~\citet{Kalbfleisch:Lawless:1989}. The corresponding approximate likelihood $\ell^{\textnormal{app}}(\boldsymbol\theta;\mathbf{f})$ reads
\begin{align*}
\begin{split}
\ell^{\textnormal{app}}(\boldsymbol\theta;\mathbf{f})
&=\sum_{j,k=1}^K \int_{V}^{C} \log \mu^\ast_{jk}\{ t ; \mathbf{J}(t-),\boldsymbol\theta \} \diff N_{jk}(t) \\
&\quad- \int_{V}^{C} I_j(t) \mu^\ast_{jk}\{t ; \mathbf{J}(t-), \boldsymbol\theta \} \textnormal{pr}_{U}[\eta-t ; \{\mathbf{J}(t-),t,k\}, \mathbf{f}] \diff t.
\end{split}
\end{align*}
Analogously to Equation~\eqref{eq:FinalLikelihood}, we introduce the imputed version $\ell^{\textnormal{app}}(\mathbf{Z},\boldsymbol\theta;\hat{\mathbf{f}}_n,\hat{\mathbf{g}}_n)$. Maximization of the approximate likelihood is significantly simpler than the true likelihood; computationally, the approximate case is the same as a standard multistate likelihood, but where the contribution of a given path is weighted with a scalar and where the exposure is reduced based on the closeness to the time of analysis via the reporting delay distribution. Practical aspects of the implementation are discussed in Section D of the Supplementary material. 
\begin{remark}
   The framework could be extended to allow for time-dependent covariates by using the associated partial likelihood. Note, however, that the covariate paths then have to be included in $\mathcal{H}$, meaning that the adjudication model also needs to depend on these covariates. Since the adjudication model has to be used prospectively, one would need a model for how the covariates evolve over time. This is about as complicated as including the covariates in the marked point process itself, which is also more natural when the goal is prediction rather than inference. 
\end{remark}

\section{Asymptotic Properties} \label{sec:Asymptotics}

The large sample properties of $\hat{\boldsymbol\theta}_n$ largely follow from the two-step M- and Z-estimation results of~\citet{Newey:McFadden:1994},~\citet{Hahn:1996}, and~\citet{Delsol:Van:2020}. We also show weak consistency and asymptotic normality of a bootstrap procedure analogous to  Efron's simple nonparametric bootstrap as studied in~\citet{Gross:Lai:1996}, where $(\mathbf{Z}^{\textnormal{obs}}_{ni})_{i=1}^n$ are drawn with replacement from $(\mathbf{Z}^{\textnormal{obs}}_{i})_{i=1}^n$. The estimators based on the bootstrap sample are denoted $\hat{\mathbf{g}}^\textnormal{boot}_n$, $\hat{\mathbf{f}}^\textnormal{boot}_n$, and $\hat{\boldsymbol\theta}^\textnormal{boot}_n$. The main result is:
\begin{theorem} \label{thm:asymp} \leavevmode
\begin{enumerate}[label=(\roman*)]
    \item Consistency: Under Assumptions 1-3 and Assumptions 5-7 from Section E of the Supplementary material, $\hat{\mathbf{g}}_n$, $\hat{\mathbf{f}}_n$, and $\hat{\boldsymbol\theta}_n$ as well as $\hat{\mathbf{g}}^\textnormal{boot}_n$, $\hat{\mathbf{f}}^\textnormal{boot}_n$, and $\hat{\boldsymbol\theta}^\textnormal{boot}_n$ are weakly consistent.
    \item Asymptotic normality: Under Assumptions 1-3 and Assumptions 5-8 from Section E of the Supplementary material, $\hat{\mathbf{g}}_n$, $\hat{\mathbf{f}}_n$, and $\hat{\boldsymbol\theta}_n$ are asymptotically normal. Furthermore $n^{1/2}(\hat{\boldsymbol\theta}_n-\boldsymbol\theta_0 ) \rightarrow N ( 0, \mathbf{V} )$ and $n^{1/2}(\hat{\boldsymbol\theta}^\textnormal{boot}_n-\hat{\boldsymbol\theta}_n ) \rightarrow N ( 0, \mathbf{V} )$ in distribution.
\end{enumerate}
\end{theorem}
\noindent The proof and the form of $\mathbf{V}$ are given in Section E of the Supplementary material. Implementing confidence intervals via estimators for $\mathbf{V}$ is challenging since there is no closed form for the score. However, Theorem~\ref{thm:asymp} implies that a percentile confidence interval based on a simple bootstrap procedure is valid. One may also use an empirical bootstrap variance estimator for $\mathbf{V}$ if $n^{1/2}(\hat{\boldsymbol\theta}^\textnormal{boot}_n-\hat{\boldsymbol\theta}_n )$ is uniformly integrable. 

\begin{remark}
    As may be seen from the proof of Theorem~\ref{thm:asymp}, the proposed estimator can, under smoothness conditions, be cast as a generalized method of moment estimator with weighting matrix equal to the identity by using score equations as moment conditions. In the class of regular asymptotically linear estimators that only restrict the statistical model through moment conditions, it is known from~\citet{Hansen:1982},~\citet{Chamberlain:1987}, and Section~5 of~\citet{Newey:McFadden:1994} that generalized method of moment estimators satisfying Equation~(5.4) of~\citet{Newey:McFadden:1994} are efficient. A sufficient condition for this is that the weighting matrix is equal to the inverse of the variance of the score, while for the maximum likelihood estimator, the identity weight matrix is also efficient due to Bartlett's identities. Despite using score-type equations as moment conditions, our proposed estimator is not efficient, since an identity of Barlett type is not generally satisfied. \citet{Hwang:Sun:2018}, however, find that the finite sample performance of the efficient estimator is often inferior to an a priori chosen weighting matrix, since the uncertainty in estimating the weighting matrix is ignored in the asymptotics. If the cost-benefit trade-off is unclear, they recommend sticking to an a priori chosen weighting matrix, in line with our proposal of using $\hat{\boldsymbol\theta}_n$. Compared to the efficient estimator, our proposal also has the computational advantage that the score need not be computed.
\end{remark}

\section{Numerical Study} \label{sec:Numerical}

\subsection{Data-generating process}

We use the multistate event model and three stage adjudication model depicted in Figure~\ref{fig:Numerical}. A jump is confirmed if and when state $3$ of the adjudication model is hit. Multistate models are simulated using Lewis' thinning algorithm from~\citet{Ogata:1981}. We refer to Section F of the Supplementary material for additional figures, tables, and details regarding the numerical study.  We compare our approach with oracle methods that use the filtration $\mathcal{F}^\ast$, which is unavailable to the statistician, and naive methods that either use all reported events without back-censoring (Naive 1) or back-censor by one year and delete reported events if they have been under adjudication for over two years at time $\eta$ (Naive 2). In Section F of the Supplementary material, we additionally compare with an estimator based on the logistic regression approach of~\citet{Cook:Kosorok:2004}; as expected, this estimator performs poorly as it is being used in a situation it was not designed to handle, confer with Section~\ref{subsec:literature}.

\begin{figure}[ht!] 
\begin{minipage}[b]{.49\linewidth}
\begin{figure}[H]
\centering
\scalebox{0.8}{
   \begin{tikzpicture}[node distance=8em, auto]
	\node[punkts] (g) {$1$};
	\node[right=4cm of g] (i1) {};
        \node[punkts, above=0.5cm of i1] (i2) {$2$};
        \node[punkts, below=0.5cm of i1] (i3) {$3$};
	\path (g) edge [pil] node [above=0.25cm]  {$\mu_{12}^\ast$} (i2)
	;
        \path (g) edge [pil] node [below=0.25cm]  {$\mu_{13}^\ast$, $U$} (i3)
	;
        \path (i2) edge [pil] node [right=0.25cm]  {$\mu_{23}^\ast$, $U$, $\xi$} (i3)
	;
    \end{tikzpicture}
}
\end{figure}
\end{minipage}
\begin{minipage}[b]{.49\linewidth}
\begin{figure}[H]
\centering
\scalebox{0.8}{
   \begin{tikzpicture}[node distance=8em, auto]
	\node[punkts] (g) {$1$};
	\node[punkts, right=2cm of g] (i1) {$2$};
        \node[punkts, right=2cm of i1] (i2) {$3$};
        \path (g) edge [pil] node [below=0.25cm]  {$\omega_{12}$} (i1)
	;
        \path (i1) edge [pil, bend left=0.15] node [below=0.25cm]  {$\omega_{23}$} (i2)
	;
    \end{tikzpicture}
}
\end{figure}
\end{minipage}
    \caption{Event history model (left) and adjudication model (right). Symbols $U$ and $\xi$ indicate the presence of reporting delays and adjudication processes, respectively.}
    \label{fig:Numerical}
\end{figure}

\noindent A total of $400$ samples of size $n=1500$ with time-horizon $\eta=5$ years are considered. The time origin is set to $0$ for all subjects, so time represents calendar time rather than time since study entry.  For a given subject, the data is generated as follows:\ $V \sim \textnormal{Uniform}(0,1)$ and $C \mid V \sim \textnormal{Uniform}(V,\eta)$. Subjects enter in state $1$ with no additional left-truncation information. The baseline covariate is $X \sim \textnormal{Uniform}(-4,4)$. We let $\boldsymbol\theta = (\log 0.15,0.1,0.4,\log 0.1, 0.03, -0.3, -0.3)$ and generate events using the rates $\Tilde{\mu}_{12}(t;X)=\exp\{\theta_{1}+\theta_{2} \times (t+X)+\theta_{3} \times\textnormal{sin}(0.5\pi X) \}$, $\Tilde{\mu}_{13}(t;X)=\exp\{\theta_{4}+\theta_{5} \times t^2+\theta_{6}\times\textnormal{cos}(0.5\pi X)\}$, and $\Tilde{\mu}_{23}(t; \Tilde{T}_1, X)=\exp\{\theta_{7} \times (t-\Tilde{T}_1) X^2\}$, which upon deletion of unreported jumps gives $(\Tilde{T}_m,\Tilde{Y}_m)_{m \geq 1}$ or upon deletion of rejected jumps gives $(T^\ast_m,Y^\ast_m)_{m \geq 1}$. The rates are relatively large, and on average around 260 transitions from 1 to 3, 415 transitions from 1 to 2, and 180 transitions from 2 to 3 are generated. 

To model the reporting delays, we assume a Weibull distribution with proportional reverse time hazards, that is $\alpha(t;X) = \alpha_0(t)\exp(X \beta)$ for $\alpha_0(t)=\frac{k \lambda^k t^{k-1}}{\exp\{(\lambda t)^k\}-1}$. This results in a power model for the distribution function, becoming $t \mapsto [1-\exp\{-(\lambda t)^k\}]^{\exp(X \beta)}$ which we denote $\textnormal{Weibull}(\lambda,k,\beta)$. Let $\mathbf{f}=(2,0.5,0.1,1,1.5,0.2)$. Reporting delays have $\textnormal{Weibull}(f_{1},f_{2},f_{3})$ distribution when coming from state $1$ and $\textnormal{Weibull}(f_{4},f_{5},f_{6})$ when coming from state $2$. The generated reporting delays have means close to one and standard deviations around $0.15$ and $0.05$, respectively. The adjudication rates are $\omega_{12}(t;\mathcal{H}_{t-}) = g_{1} \times \{X/(t-\Tilde{T}_2-\Tilde{U}_2 +2)\}^2$ and $\omega_{23}(t;\mathcal{H}_{t-}) = \exp\{g_{2} \times (t-\tau_{1,1})\}$ with $\mathbf{g} = (0.8,-1.2)$, leading to a long-term confirmation rate of $37\%$. To go from $\Tilde{\mu}_{jk}$ to $\mu^\ast_{jk}$, note that $\mu^\ast_{12}(t;X)=\Tilde{\mu}_{12}(t;X)$ and $\mu^\ast_{13}(t;X)=\Tilde{\mu}_{13}(t;X)$, while an expression for $\mu^\ast_{23}\{t,D(t);X\}$ may be found in Section F of the Supplementary material.

\subsection{Finite sample performance}

The regularity assumptions from Theorem~\ref{thm:asymp} are easily seen to be satisfied, and we hence expect the proposed estimators to be unbiased and asymptotically normal. In Table~\ref{table:theta} we report the bias and the empirical standard deviation (SD) of the parameter estimators of $\boldsymbol\theta$ using the different methods. A corresponding table for the estimators $\hat{\mathbf{g}}_n$ and $\hat{\mathbf{f}}_n$ can be found in Section F of the Supplementary material. The results are generally as expected. We find comparable bias for the proposed method and oracle methods, with a higher SD for the proposed method. The bias is slightly higher for the Poisson approximation with SD comparable to the proposed method. The naive methods are slightly worse than the proposed method for $(\theta_1,\theta_2,\theta_3)$, moderately worse for $(\theta_4,\theta_5,\theta_6)$, and substantially worse for $\theta_7$. The naive methods generally have comparable performance except that Naive 1 is superior for $\theta_4$ and Naive 2 is superior for $\theta_7$. We conclude that our proposed method works well. 

The Poisson approximation is also expected to perform well since there is no jump after the transition affected by reporting delays, confer with Section~\ref{sec:est_theta}. Since the hazards are relatively large, the decrease in bias when going from the Poisson approximation to our proposed method is non-negligible. Depending on the application, this reduction in bias might be worth the significantly longer computation times. With our implementation and hardware setup, the Poisson approximation for $\theta_\ell$ $(\ell=1,\dots,6)$ took only 7 seconds, while the proposed method took about 300 seconds.

We also seek to verify the bootstrap procedure outlined in Section~\ref{sec:Asymptotics}. We here focus on $\theta_7$ since its estimation is affected by both reporting delays and incomplete event adjudication. Based on $k=400$ estimates of $\theta_7$ with $1000$ bootstrap resamples and confidence level $1-\alpha=0.90,0.95,0.99$, we obtain coverage rates of $89.5, 95.0,98.7$. We conclude that the bootstrap works well as the coverage rates are close to the nominal level.

\begin{table}
\centering
\caption{Bias and empirical standard deviation (SD) of the estimator $\hat{\boldsymbol\theta}_n$ based on 400 simulations of size $n=1500$.}
\begin{tabular}{lccccccccccccccc} 
\hline
 &  \multicolumn{2}{c}{Proposed method}   &  \multicolumn{2}{c}{Oracle}  &  \multicolumn{2}{c}{Poisson approx.}  & \multicolumn{2}{c}{Naive 1} & \multicolumn{2}{c}{Naive 2}  \\
  \cmidrule(lr){2-3} \cmidrule(lr){4-5} \cmidrule(lr){6-7} \cmidrule(lr){8-9} \cmidrule(lr){10-11}
Parameter       &  Bias   & SD    &  Bias  & SD      & Bias    & SD       & Bias   & SD      & Bias    & SD \\[3pt]
 \hline
$\theta_1=\log 0.15$ & -.004   & .067  &  -.005 & \;.031    & -.010   & \;.067     & -.011  & \;.066    & -.010   & \;.066   \\
$\theta_2=0.1$ & -.000   & .020  & -.001  & \;.020    & -.006   & \;.020     & -.006  & \;.020    & -.006   & \;.021  \\
$\theta_3=0.4$ & \;.003    & \;.078  & \;.003   & \;.078   & -.002   & \;.078     & -.002  & \;.078    & -.000   & \;.079  \\
$\theta_4=\log 0.1$ & \;.003    & \;.084  & \;.001   & \;.083    & \;.012    &  \;.091    & -.041  & \;.077    & -.051   & \;.082  \\
$\theta_5 = 0.03$ & \;.000    & \;.012  & -.001  & \;.013    & -.006   & .016     & -.018  & \;.011    & -.015   & \;.014  \\
$\theta_6 = -0.3$ & -.000   & \;.094  & -.001  & \;.088    & \;.007    & .094     & -.009  & \;.087    & -.007   & \;.090  \\
$\theta_7 = -0.3$ & -.011   & \;.066  &  -.011 & \;.054    & -.012   & .066     & \;.157   & \;.023    & \;.148    & \;.069 \\
 \hline
\end{tabular}
\label{table:theta}
\end{table}

\subsection{Predictive accuracy and robustness against misspecification}

To evaluate the predictive accuracy beyond individual parameter estimates, we explore the performance of the different approaches for a relevant estimand. The estimand of choice is the expected duration spent in state $2$ before time $t$, which we denote $V_a(t;X) = \mathbb{E}\{ \int_{[0,t] } I^\ast_2(s) \diff s \mid X \}$. In an actuarial context, this may be recognized as an active reserve for a unit disability annuity if states (1), (2), and (3) are identified with the active, disabled, and dead states, respectively; in the notation of Section~\ref{sec:Motivation}, it corresponds to $B(\diff t) = I^\ast_2(t) \diff t$. We use the same simulations as before and plug in the estimated parameter values. The estimand $V_a$ is calculated by solving the system of differential equations from Corollary~7.2 in~\citet{Adekambi:Christiansen:2017} using the fourth-order Runge-Kutta method, see also Section 8 in their paper for a discussion of numerical solution methods. For a given estimator $\hat{V}_a$, we compute the error metrics 
\begin{align*}
    \textnormal{MSE}(\hat{V}_a,V_a) &= \int_{[0,\eta]} \int_{[-4,4]} \{ \hat{V}_a(s;x) - V_a(s;x)\}^2 \diff s \diff x, \\
    \textnormal{MAE}(\hat{V}_a,V_a) &= \int_{[0,\eta]} \int_{[-4,4]} \vert \hat{V}_a(s;x) - V_a(s;x) \vert \diff s \diff x
\end{align*}
in each simulation and then average the results. The resulting values are reported in Table~\ref{table:prederr} and are as expected. Similar to Table~\ref{table:theta}, we see that the proposed method performs just a tad worse than the oracle method, that the Poisson approximation performs slightly worse, and that the naive methods perform much worse. 
\begin{table}
\centering
\caption{Average MSE and MAE based on 400 simulations of size $n=1500$.}
\begin{tabular}{ccccccccccccccc} 
\hline
 \multicolumn{2}{c}{Proposed method}   &  \multicolumn{2}{c}{Oracle}  &  \multicolumn{2}{c}{Poisson approx.}  & \multicolumn{2}{c}{Naive 1} & \multicolumn{2}{c}{Naive 2}  \\
  \cmidrule(lr){1-2} \cmidrule(lr){3-4} \cmidrule(lr){5-6} \cmidrule(lr){7-8} \cmidrule(lr){9-10}
        MSE   & MAE    &  MSE  & MAE      & MSE    & MAE       & MSE   & MAE      & MSE    & MAE \\[3pt]
 \hline
.084   & 1.18  &  .083 & 1.17    & .089   & 1.21     & .149  & 1.53    & .151   & 1.55   \\
 \hline
\end{tabular}
\label{table:prederr}
\end{table}

We also investigate the performance of the proposed method under misspecification. Here the parameters are estimated under a misspecified parametric family for $(T_m^\ast,Y_m^\ast)_{m \geq 1}$ characterized by the hazard rates
\begin{align*}
    \mu^{\textnormal{miss}}_{12}(t;X)&=\exp\{\theta_{1}+\theta_{2} \times t+\theta_{3} \times X \}, \\
    \mu_{13}^{\textnormal{miss}}(t;X)&=\exp\{\theta_{4}+\theta_{5} \times t+\theta_{6}\times X\}, \\
    \mu_{23}^{\textnormal{miss}}(t; \Tilde{T}_1, X)&=\exp\{\theta_{7} +  \theta_{8} \times (t-\Tilde{T}_1) + \theta_{9} \times X\}.
\end{align*}
This misspecification is substantial since it disregards the non-linearities and interactions present in the true hazards. We both consider the case where the reporting delay and adjudication models are correctly specified, and a doubly misspecified case where these are also misspecified via the following parametric families:
\begin{align*}
    \omega_{12}^{\textnormal{miss}}(t;\mathcal{H}_{t-}) &= \exp\{ g_1\times X+g_2 \times (t-\tilde{T}_2-\tilde{U}_2) \}, \\
    \omega_{23}^{\textnormal{miss}}(t;\mathcal{H}_{t-}) &= \exp \{ g_3\times X+g_4 \times (t-\tau_{1,1}) \},
\end{align*}
and $\alpha(t;X)=\alpha_0(t)\exp(X\beta)$ with $\alpha_0$ being the reverse time hazard for a Gamma distribution with shape parameter $k$ and rate parameter $\lambda$.

\begin{figure}
    \centering
    \includegraphics[scale=0.65]{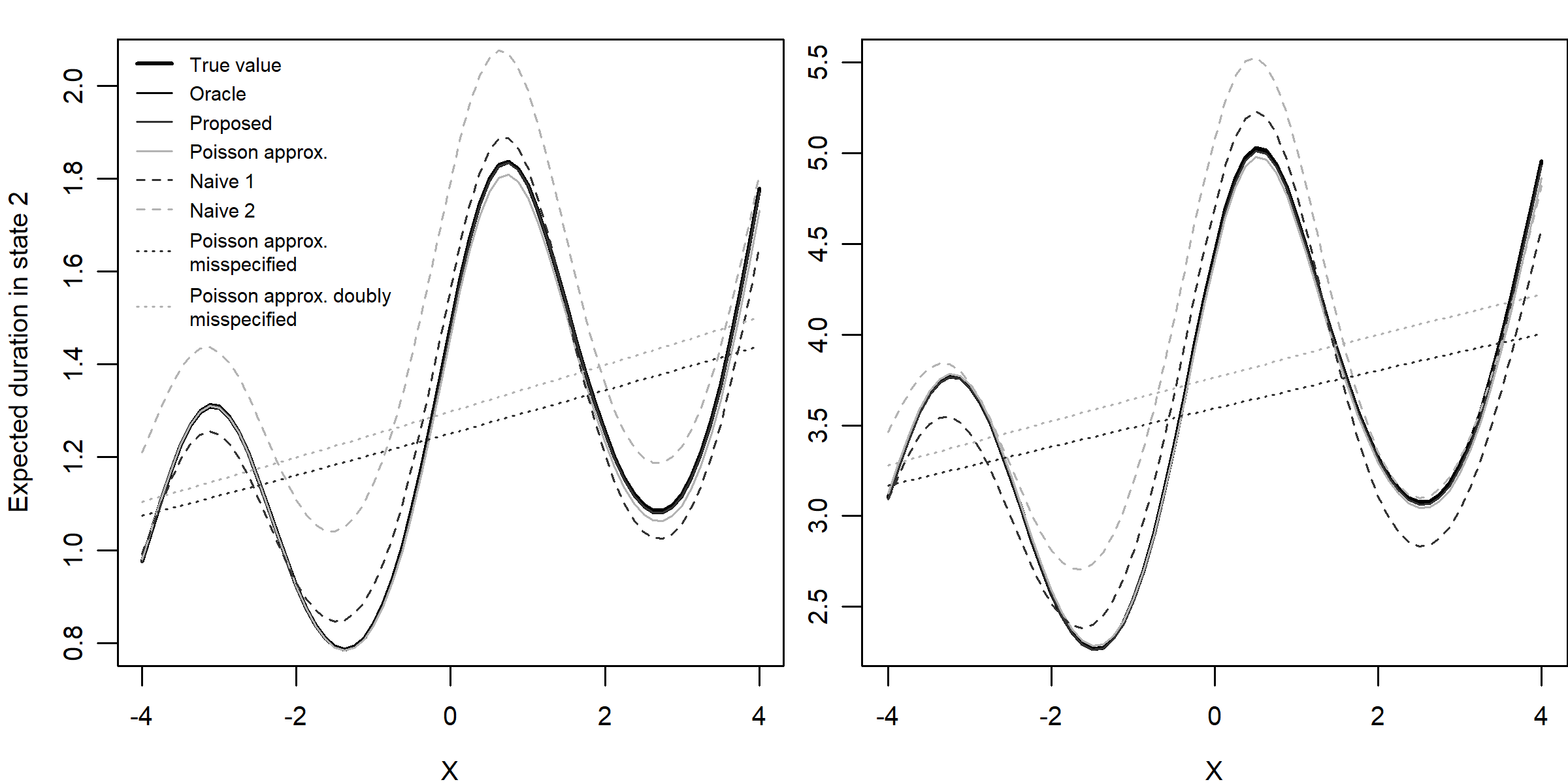}
    \caption{Plot of plug-in estimates of $V_a(t;x)$ averaged over 400 samples of size $n=1500$ with $t=\eta$ on the left and $t=2\eta$ on the right.}
    \label{fig:VaPlot}
\end{figure}

It is not meaningful to compare the parameter values for the misspecified and correctly specified hazards, so we solely focus on the predictive accuracy for $V_a(t;X)$ at time $t=\eta$ and $t=2\eta$. Note that the case $t=2\eta$ requires extrapolation, which is not inherently an issue for parametric models. We observed that the Poisson approximation performed very well for this estimand, as can be seen in Figure~\ref{fig:VaPlot}, so we have chosen to only consider the misspecified model used in conjunction with the Poisson approximation. The results are illustrated in Figure~\ref{fig:VaPlot}.

The oracle, proposed, and Poisson approximation approaches are all very close to the true value of $V_a$ leading their graphs to overlap almost perfectly, with the Poisson approximation showing slight deviations around $X=0$ and $X=2$. Naive 1 performs reasonably albeit with some noticeable bias, while Naive 2 seems to be systematically biased upwards except for $t=2\eta$ and $X \geq 2$ where it is very close to the true value. The misspecified parametric model captures the local $X$-dependence of $V_a(t;X)$ the worst, but identifies the correct overall level for the estimand. This behavior is similar to how maximum likelihood estimators behave under misspecification, see~\citet{Halbert:1982}. Including polynomial or sinusoidal transformations of $X$ in the parametric model would likely have led to a better fit for the local behavior of $V_a$, but including too many non-linear terms could lead to high variance estimators due to overfitting. The performance is comparable for $t=\eta$ and $t=2\eta$, showing that our method is able to capture the general trends reasonably well even under substantial model misspecification.

\newpage
\section{Data Application} \label{sec:Application}

\subsection{Data characteristics}
We introduce a new data set LEC-DK19 (Loss of Earning Capacity -- Denmark 2019) collected by a large Danish insurer in the period 31/01/2015 to 01/09/2019 with time of analysis $\eta=\textnormal{01/09/2019}$ and apply our proposed estimation procedure to this data. The data set is closely related to the insurance problem discussed in Section~\ref{sec:Motivation}. Observations are available on a monthly grid. Subjects first appear in the data when they enter the portfolio and are censored when they leave the portfolio. Claims may be reported after censoring, meaning that information on reporting delays and adjudication may arrive after censoring. A total of 416,483 insured are included across five tables concerning disabilities, reactivations, disability delays, disability adjudications, and reactivation adjudications. Of the 413,139 insured in the disability data, 1,773 (0.43\%) have an unrefuted disability; disabilities are thus rare. Of the 3,011 insured in the reactivation data, 1,133 (37.63\%) have an unrefuted reactivation. The data is based on raw data that has been anonymized and slightly altered so as not to reveal any confidential information about the individual subjects or the insurance portfolio. Available baseline covariates are gender and date of birth.

\begin{figure}[ht!]
\begin{minipage}[b]{.49\linewidth}
\begin{figure}[H]
\centering
\scalebox{0.8}{
   \begin{tikzpicture}[node distance=8em, auto]
	\node[punkts] (g) {$1$};
	\node[punkts, right=2cm of g] (i1) {$2$};
        \node[punkts, right=2cm of i1] (i2) {$3$};
        \node[punkts, below=1.25cm of i1] (d) {$4$};
        \path (g) edge [pil] node [above=0.2cm]  {$\mu^\ast_{12},U,\xi$} (i1)
	;
        \path (i1) edge [pil, bend left=0.15] node [above=0.2cm]  {$\mu^\ast_{23},\xi$} (i2)
	;
          \path (g) edge [pil, bend left=0.15] node [below=0.25cm]  {$\mu^\ast_{14}$} (d)
	;
         \path (i1) edge [pil, bend left=0.15] node [left=0.15cm]  {$\mu^\ast_{24}$} (d)
	;
         \path (i2) edge [pil, bend left=0.15] node [below=0.25cm]  {$\mu^\ast_{34}$} (d)
	;
    \end{tikzpicture}
}
\end{figure}
\end{minipage}
\begin{minipage}[b]{.49\linewidth}
\begin{figure}[H]
\centering
\scalebox{0.8}{
   \begin{tikzpicture}[node distance=8em, auto]
	\node[punkts] (g) {$3$};
	\node[punkts, right=2cm of g] (i1) {$1$};
        \node[punkts, right=2cm of i1] (i2) {$2$};
        \node[, right=0.8cm of i1] (a) {};
        \node[punkts, below=1.5cm of a] (a2) {$4$};
         \path (i1) edge [pil] node [above=0.155cm]  {$\omega_{13}$} (g)
	;
        \path (i1) edge [pil, bend right=20] node [below=0.15cm]  {$\omega_{12}$} (i2)
	;
         \path (i2) edge [pil, bend right=20] node [above=0.1cm]  {$\omega_{21}$} (i1)
	;
         \path (i1) edge [pil] node [left=0.15cm]  {$\omega_{14}$} (a2)
	;
         \path (i2) edge [pil] node [right=0.15cm]  {$\omega_{24}$} (a2)
	;
    \end{tikzpicture}
}
\end{figure}
\end{minipage}
    \caption{Event history model (left) and adjudication model (right). For events, active is 1, disabled is 2, reactivated is 3, and dead is 4. For adjudications, active report is 1, inactive report is 2, adjudicated is 3, and dead is 4.}
    \label{fig:realdata}
\end{figure}

\subsection{Model specification}
We model the data with the multistate model depicted in Figure~\ref{fig:realdata}. Subjects can become disabled, reactivate from disablement, and die. A disability starts in the adjudication state 1, is confirmed if state 3 is hit, and otherwise rejected. A reactivation starts in adjudication state 2 and the reactivation is annulled if state 3 is hit. Only disablements seem to exhibit reporting delay, likely because the reactivations result from the insurer terminating payments. 

We employ the Poisson approximation, which implies that the hazards can be estimated separately. We here only estimate $\mu^\ast_{12}$ and $\mu^\ast_{23}$, partly because these are of primary interest, but also because only deaths recorded during the adjudication period are included in the data, impeding estimation of the full event history model. Via back-censoring, it is seen that the hazard $\mu^\ast_{12}$ likely has order of magnitude $10^{-2}$, which is one to two orders of magnitudes less than the hazards in the numerical study of Section~\ref{sec:Numerical}, so the incurred approximation error for $\mu^\ast_{12}$ is expected to be negligible. The lack of mortality data makes it difficult to assess the approximation error for $\mu^\ast_{23}$, however an alternative is to approximate the disabled mortality hazard with zero. 

The hazards are regressed on the baseline covariates and the time elapsed since 31/01/2015, and $\mu^\ast_{23}$ is additionally regressed on the duration spent in the disabled state. An adjudication hazard is estimated for each transition. Following the form of $\mathcal{H}_t$ suggested in Section~\ref{sec:Model}, the adjudication hazards for disabilities are regressed on the duration since the disability occurred and the associated reporting delay, and the adjudication hazards for reactivations are regressed on the time since the disability and reactivation events. In addition, for disabilities we regress on whether state 2 has previously been hit. For each adjudication and event hazard, the covariates enter in a linear predictor with log link. The log link implies that the hazards are always positive, so a disability rejection or reactivation confirmation is only certain once the insured dies, but the probability of a confirmed jump decreases towards zero as the duration of adjudication increases towards infinity. For reporting delays, we employ the Weibull proportional reverse time hazards of Section~\ref{sec:Numerical} with age entering as age at disability onset. Since the data set contains monthly records, we employ the Poisson regression approach for implementation as described in Section D of the Supplementary material.

\subsection{Empirical results}
In Table~\ref{table:dataTheta}, we present the non-nuisance parameter estimates and percentile bootstrap confidence intervals computed using 400 bootstrap samples. For comparison, we additionally employ a naive estimation procedure consisting of using all reported unrefuted events with no back-censoring  and with one year of back-censoring. A corresponding table for the nuisance parameters is provided in Section G of the Supplementary material. For this data, we see that adjusting for reporting delays and incomplete event adjudication has a substantial effect on the estimated level and calendar-time dependence of the disability hazard. The effects on the parameter estimates of the reactivation hazard are also noticeable, but contrary to the case for the disability hazard, the naive estimates stay firmly within the confidence intervals except for the calendar-time dependence. To further illustrate the results, we plot transition probabilities based on the estimated hazards in Section G of the Supplementary material.  

In Figure~\ref{fig:FittedRates}, predicted rates across the data set, using both the proposed and naive procedures, are shown against the corresponding occurrence-exposure rates aggregated by yearly tenths. Despite the smooth form of the conditional hazards, the aggregation with respect to covariates produces jagged curves as the mix of covariates changes over time. Since the naive occurrence-exposure rates use all unrefuted events, they are higher than the adjusted occurrence-exposure rates for reactivations. For disabilities, they are higher in the first years and lower in the later years, where the influence of reporting delays is more substantial. For the disability rate, back-censoring by one year reduces the effect of reporting delays substantially. Nevertheless, the predictions are noticeably lower than the proposed method since the adjusted occurrence-exposure rates indicate an even larger increase in the final year of the observation window. For the reactivation rate, back-censoring by one year leads to a higher rate than not back-censoring. This is likely because reactivation adjudications often are longer than one year, so part of the effect remains, and because the reactivation rate seems to be substantially lower in the final year.
Both rates show large calendar-time dependence, which is a distinct feature of the data set that implies an especially poor forecasting performance of back-censoring and underscores the importance of our proposed methods for disability insurance applications.

\begin{table}
\centering
\caption{Parameter estimates (Est.), Naive parameter estimates (Naive), Back-censoring parameter estimates (Backcens.), and 95\% bootstrap percentile confidence interval (CI) for  $\boldsymbol\theta$ using the proposed method and 400 bootstrap resamples.}
\begin{tabular}{lcccccccc}
\hline
 &  \multicolumn{4}{c}{$\mu^\ast_{12}$} &  \multicolumn{4}{c}{$\mu^\ast_{23}$} \\
  \cmidrule(lr){2-5} \cmidrule(lr){6-9}
Parameter &  Est.  & Naive & Backcens. & CI &  Est. & Naive & Backcens. & CI \\[3pt]
 \hline
Age & .023  & .022 & .024 & (0.018, 0.027) &  -.012  & -.012 & -.012 & (-0.018, -0.007)    \\
Male & -7.46  & -7.11 & -7.29 & (-7.65, -7.23) & .334 & .262 & .290 & (0.073, 0.624) \\
Female & -8.66 &  -8.25 & -8.78 & (-8.87, -8.40) & .756  & .674 & .673 & (0.456, 1.07) \\
Time & .304 & .098 & .198 & (0.188, 0.357) & -.125 & -.068 & -.046 & (-0.168, -0.073) \\
Duration & -  & - & - & -  & -1.04  & -1.03 & -1.10  & (-1.15, -0.942) \\
 \hline
\end{tabular}
\label{table:dataTheta}
\end{table}

\begin{figure}
    \centering
    \includegraphics[scale=0.56]{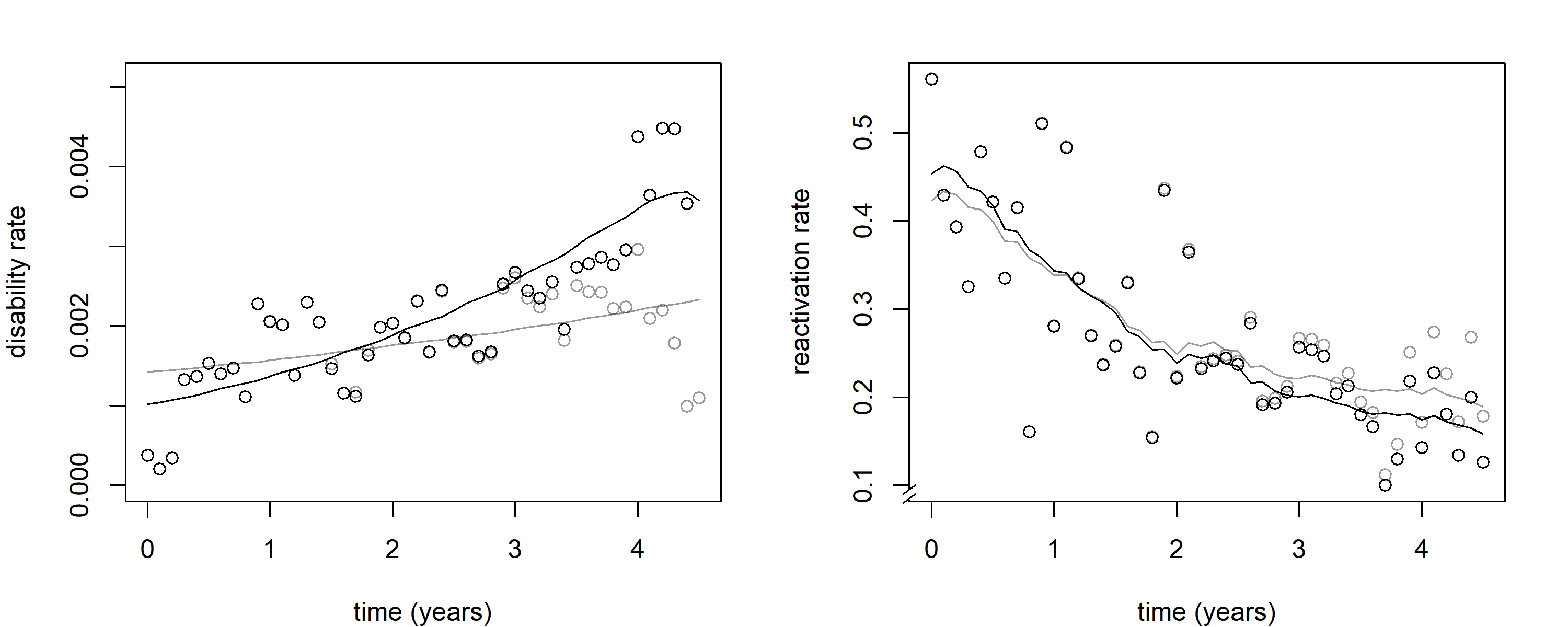}
    \caption{Fitted rates (lines) and occurrence-exposure rates (points) for the proposed method (black), the naive method with no back-censoring (light gray), and one year of back-censoring (dark gray). The dark gray line is dashed in the final year of the observation window to signify that it involves extrapolation.  Disability rates are shown on the left and reactivation rates on the right.}
    \label{fig:FittedRates}
\end{figure}

 The alternative approach, where the Poisson approximation is not employed but the disabled mortality is approximated with zero, leads to the estimated regressors of $\mu_{23}^\ast$ being $(-0.011,0.373,0.840,-0.260,-0.907)$. Compared with the regressors obtained when using the Poisson approximation, the main difference is that a portion of the duration dependence is shifted to calendar-time dependence. This would have a noticeable influence on the forecasted reactivation rates. For real applications, we therefore recommend using the full estimation procedure with access to mortality data when estimating the reactivation hazard. We also recommend validating the models by comparing the original predictions with predictions obtained using the same data but with different time of analyses. However, this is not possible with the current data, as the disability and reactivation occurrences, along with exposures, are available only for a single valuation date. Moreover, their values at different valuation dates cannot be inferred from the adjudication data, since the same ID does not correspond to the same insured individual across the separate data tables. 
 
 Note also that usual out-of-sample validation is not applicable since the data itself is biased. Nevertheless, cross-validation based on the adjusted occurrences and exposures corresponding to the black points in Figure~\ref{fig:FittedRates} (constructed according to Section D of the Supplementary material) should still provide a valid way to choose between competing models for $\mu^\ast_{jk}$. This could be explored in the current analysis by varying the parametric form for each hazard.

\begin{figure}
    \centering
    \includegraphics[scale=0.45]{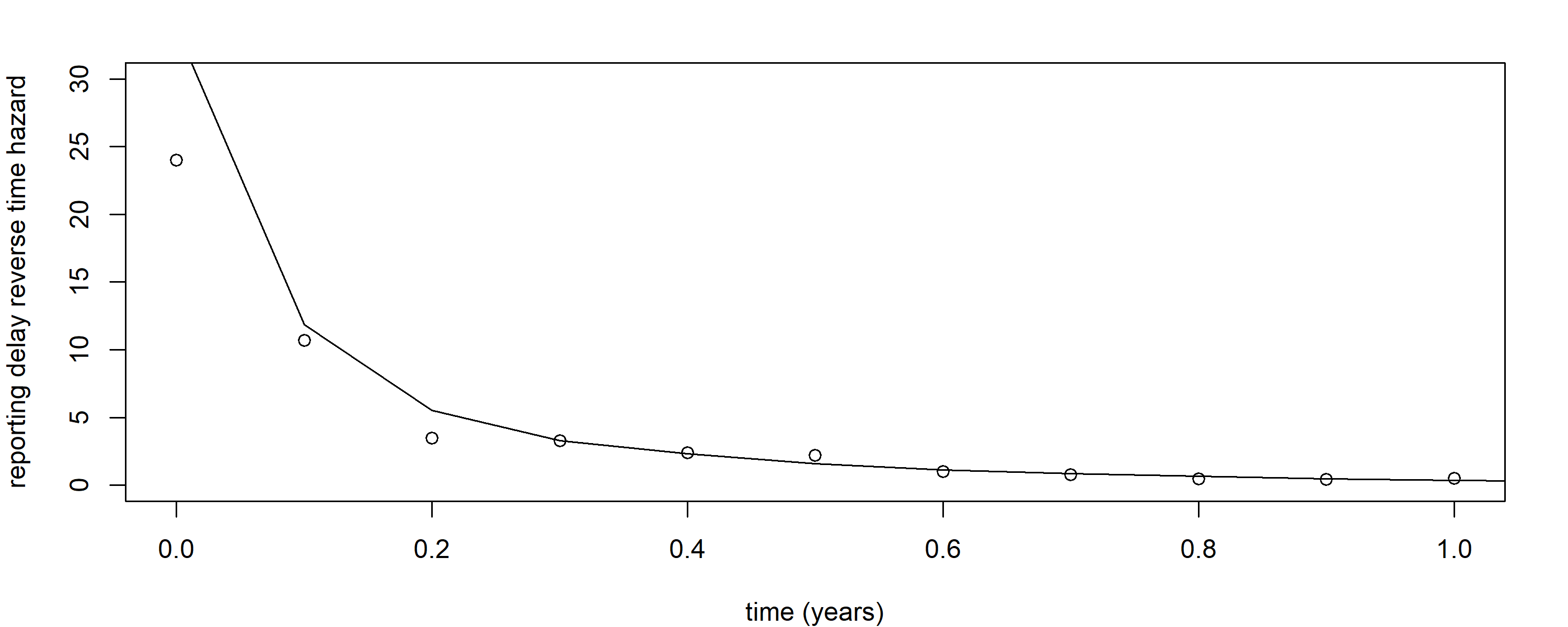}
    \caption{Fitted rate (line) and occurrence-exposure rate (points) for the reverse time hazard of the reporting delay distribution.}
    \label{fig:DelayRates}
\end{figure}

\begin{figure}
    \centering
    \includegraphics[scale=0.725]{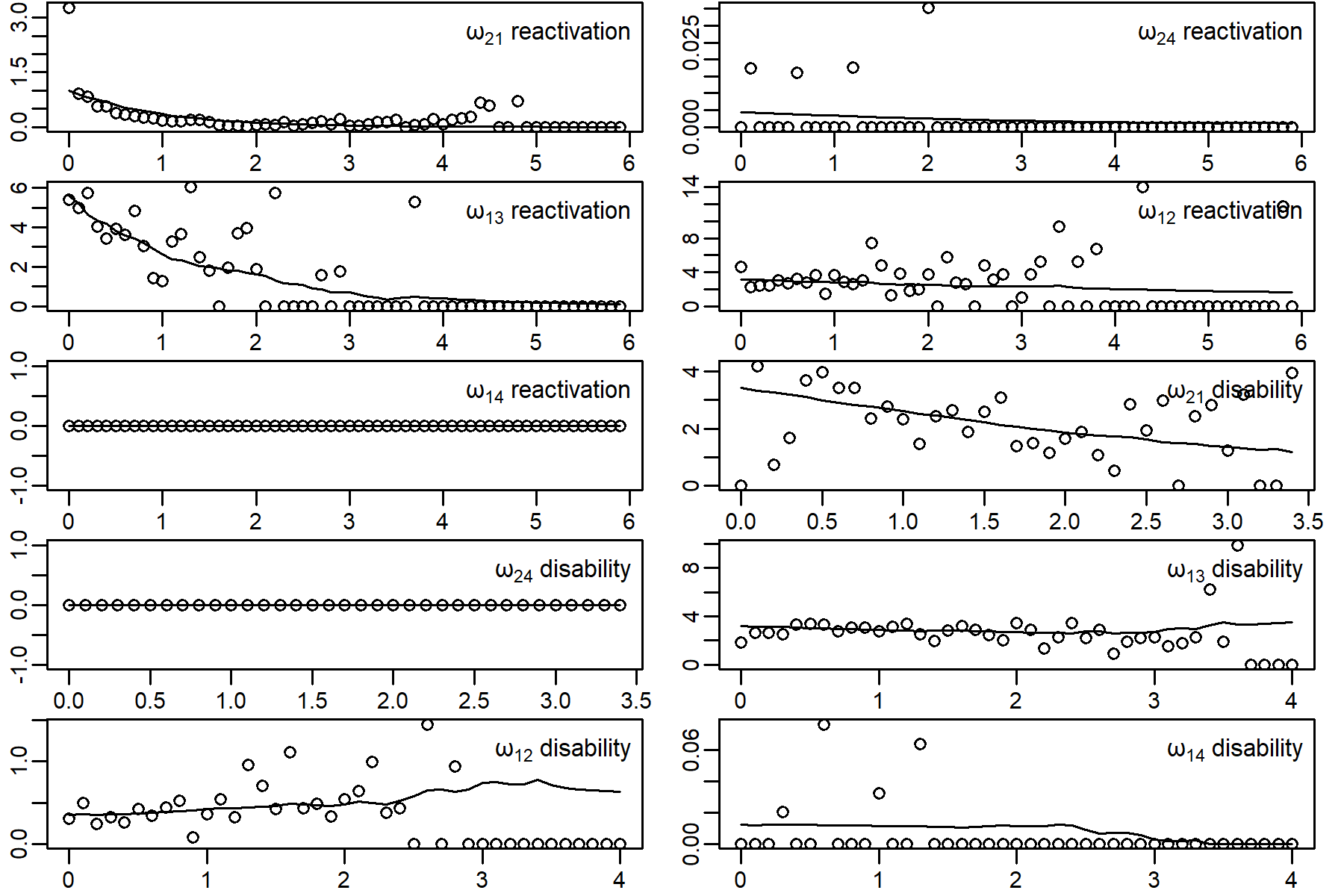}
    \caption{Fitted rates (lines) and occurrence-exposure rates (points) for the adjudication hazards graphed as a function of the duration since the unadjudicated event occurred.}
    \label{fig:AdjRates}
\end{figure}

In Figure~\ref{fig:DelayRates}, we use the link between reverse time hazard estimation and Poisson regression, which is described in Section D of the Supplementary material, to plot the fitted reverse time hazard rates against empirical occurrence-exposure rates. The model shows no obvious lack of goodness-of-fit. Figure~\ref{fig:AdjRates} contains similar plots for the adjudication hazards as a function of the duration since the unadjudicated event occurred and which also do not show no systematic deviations.

\section{Closing Remarks} \label{sec:Closing}
We have proposed a parametric two-step method to estimate hazards when data is left-truncated, right-censored, and contaminated by reporting delays and incomplete event adjudication. We proved that the estimators and their bootstrapped versions are weakly consistent and asymptotically normal. The numerical study showed favorable performance of the proposed method compared to other alternatives, but also highlighted that for small hazards the Poisson approximation performs reasonably while being considerably less computationally demanding. Our approach overcomes the need for back-censoring, meaning that the most recent data may be used. This feature is particularly useful in monitoring the emergence of new trends on a population level, confer with the data application in Section~\ref{sec:Application}.

Section A–G of the Supplementary material, Tables, and Figures referenced in Sections~\ref{sec:Model}-\ref{sec:Application} are available in Supplement A of the Supplementary material (\citet{Buchardt.etal:2025}). The R code and the LEC-DK19 data are available in Supplement B of the Supplementary Material and on GitHub: \url{https://github.com/oliversandqvist/Web-appendix-estimation-contamination}.


\begin{acks}[Acknowledgments]
  Oliver Lunding Sandqvist is an Industrial PhD student at PFA Pension and Kristian Buchardt is an Affiliated Professor of Actuarial Mathematics at the University of Copenhagen. Significant parts of the research presented in this paper were carried out while Kristian Buchardt was employed at PFA Pension. The authors thank Munir Hiabu for insightful comments that helped improve the manuscript. \vspace*{-8pt}
\end{acks}


\begin{funding}
Oliver Lunding Sandqvist's research has partly been funded by the Innovation Fund Denmark (IFD) under File No.\ 1044-00144B.
\end{funding}


\begin{supplement}
\stitle{Supplement A: Proofs and additional details.}
\sdescription{We provide proofs and additional details for the methods proposed in this paper.}
\end{supplement}

\begin{supplement}
\stitle{Supplement B: Code and data.}
\sdescription{We provide the R implementation and data used in this paper.}
\end{supplement}




\end{document}


\begin{frontmatter}
\title{Supplementary material for ‘Estimation for Multistate Models Subject to Reporting Delays and Incomplete Event Adjudication With Application to Disability Insurance’}
\runtitle{Estimation with Reporting Delays and Incomplete Event Adjudication}

\begin{aug}
\author{Kristian Buchardt},
\author{Christian Furrer}
\and
\author{Oliver Sandqvist}
\end{aug}

\begin{abstract}

The supplementary material is organized as follows. Section{\scshape~\ref{supplement:Imputation}} provides a discussion of the link between the proposed estimation procedure and imputation. In Section{\scshape~\ref{supplement:Distribution}}, we prove Lemma~\ref{lemma:hazard_link} and derive an error-bound for the Poisson approximation. In Section{\scshape~\ref{supplement:Delays}}, we discuss the estimation procedure and its efficiency under endogenous reporting delays. Section{\scshape~\ref{supplement:Implement}} provides a way to implement parts of the proposed estimation procedure via Poisson regression and explains how this link may be used to construct goodness-of-fit plots. In Section{\scshape~\ref{supplement:Asymptotics}}, we derive the large-sample behavior of our estimators and prove Theorem~\ref{thm:asymp}. Section{\scshape~\ref{supplement:Numerical}} contains further details regarding the data-generating mechanism used in the numerical study as well as additional figures and tables. Section{\scshape~\ref{supplement:Application}} contains additional figures and tables for the data application.
\end{abstract}

\end{frontmatter}


\section{Relation to Imputation} \label{supplement:Imputation}

The proposed estimation procedure for handling adjudications corresponds to frequentist imputation when the number of imputations $b$ tends to infinity. The main challenge regarding this aspect of the estimation procedure is thus that the time-dependent nature of both the events of interest and the adjudication information causes the adjudication outcomes to have a complicated distribution.

If one had instead used a finite $b$, one may form a single `clustered data' objective function $\frac{1}{b}\sum_{j=1}^b \frac{1}{n}\sum_{i=1}^n \ell(\boldsymbol\theta;\mathbf{f},k_{ij})$, where $k_{ij}$ is the imputed value of the adjudication outcome for individual $i$ in the imputed dataset $j$. Note that this does not simplify the main challenge described above, as this approach requires one to simulate the adjudication outcomes given the available information. This `clustered data' procedure is described and analyzed in~\citet{Wang:Robins:1998},~\citet{Robins:Wang:2000}, and~\citet{Kim:2011} for both finite $b$ and in the limit $b\rightarrow\infty$. An approximate maximizer of $\sum_{i=1}^n \ell(\mathbf{Z}_i,\boldsymbol\theta;\mathbf{f},\mathbf{g})$ could thus alternatively have been based on a large but finite number of imputed datasets (simulating $k_{ij}$ based on its distribution under $\mathbf{g}$ conditional on $\mathbf{Z}_i$), where the user inputs the $b \times n$ observations as `independent' observations into a standard algorithm that is able to maximize $\ell(\boldsymbol\theta;\mathbf{f})$.  For finite $b$, one may alternatively use the multiple imputation approach described on p.\ 76 of~\citet{Rubin:1987}, where an estimate $\hat{\boldsymbol\theta}_{nj}$ ($j=1,\dots,b$) is computed for each of the $b$ imputed datasets by generating the imputed values under the posterior distribution of the missing values given the observed values under some specified Bayesian model. A final estimate is then formed by $\hat{\boldsymbol\theta}_n = \frac{1}{b}\sum_{j=1}^b \hat{\boldsymbol\theta}_{nj}$. Under some regularity assumptions, p.\ 117 of~\citet{Robins:Wang:2000} shows these procedures to be $n$-asymptotically equivalent for $b=\infty$, while for finite $b$ they have the same asymptotic mean but the frequentist method has a strictly smaller asymptotic variance. 

\section{Distribution of the Thinned Marked Point Process} \label{supplement:Distribution}

\subsection{Proof of Lemma~\ref{lemma:hazard_link}}

We derive an operational expression for the hazard $\mu_{jk}$, which in turn gives an expression for the partial log-likelihood $\ell(\boldsymbol\theta;\mathbf{f})$. The calculations become straightforward by focusing on the jump-time and jump-mark distribution of the next jump in $(T_{m},Y_{m})_{m \geq 1}$ given the previous jump-times and jump-marks. Let $\mathbf{J}^\ast_m = (T^\ast_\ell,Y^\ast_\ell)_{1 \leq \ell \leq m}$ and $\mathbf{J}_m = (T_\ell,Y_\ell)_{1 \leq \ell \leq m}$ denote the first $m$ jump-times and jump-marks. The baseline covariates are suppressed throughout.

\begin{proof}
Let $P$ be the distribution under $\boldsymbol\theta$ and $\mathbf{f}$. All calculations hold on the event $(T_{m} < \infty)$. For the jump-time distribution, Assumption~\ref{AssumptionMonotoneReporting} implies that
\begin{align*}
    &P(T_{m+1}  \leq t \mid \mathbf{J}_{m} ) \\
    &=  P(T^\ast_{m+1} \leq t, T^\ast_{m+1}+U^\ast_{m+1} \leq \eta \mid \mathbf{J}^\ast_{m}, T_m^\ast+U_m^\ast \leq \eta) \\
    &= P(T_m^\ast+U_m^\ast \leq \eta \mid \mathbf{J}^\ast_{m})^{-1} \times \int_0^t \int_0^{\eta-s} P(U^\ast_{m+1} \in \diff u \mid \mathbf{J}^\ast_{m}, T^\ast_{m+1}=s) P(T^\ast_{m+1} \in \diff s \mid \mathbf{J}^\ast_{m}).
\end{align*}
For $P(U^\ast_{m+1} \in \diff u \mid \mathbf{J}^\ast_{m}, T^\ast_{m+1}=s)$, one may use the law of total probability with respect to $Y^\ast_{m+1}$ to get
\begin{align*}
    &P(U^\ast_{m+1} \in \diff u \mid \mathbf{J}^\ast_{m}, T^\ast_{m+1}=s) = \sum_{k=1}^J \textnormal{pr}_{U} \{ \diff u ; (\mathbf{J}^\ast_{m},s,k),\mathbf{f} \} \frac{\mu^\ast_{Y_{m}^\ast k}(s; \mathbf{J}_{m}^\ast, \boldsymbol\theta)}{\mu^\ast_{Y_{m}^\ast \bigcdot}(s; \mathbf{J}_{m}^\ast, \boldsymbol\theta)}.
\end{align*}
Furthermore,
\begin{align*}
    P(T^\ast_{m+1} \in \diff s \mid \mathbf{J}_{m}^\ast) &= 1_{(s \geq T^\ast_{m})} \exp\left \{-\int_{T_{m}}^s \mu^\ast_{Y_{m}^\ast \bigcdot}(v ; \mathbf{J}_{m}^\ast, \boldsymbol\theta) \diff v \right \} \mu^\ast_{Y_{m}^\ast \bigcdot }(s; \mathbf{J}_{m}^\ast,\boldsymbol\theta) \diff s
\end{align*}
by Assumption~\ref{AssumptionMonotoneReporting}. Inserting these results gives
\begin{align*}
    P(T_{m+1}  \leq t \mid \mathbf{J}_{m} ) &= \textnormal{pr}_U(\eta - T_m^\ast ; \mathbf{J}^\ast_{m}, \mathbf{f})^{-1} \times \int_{T^\ast_{m}}^t \exp\left \{-\int_{T_{m}}^s \mu^\ast_{Y_{m}^\ast \bigcdot}(v ; \mathbf{J}_{m}^\ast, \boldsymbol\theta) \diff v \right \} \\
    & \qquad \quad \times \sum_{k=1}^J \textnormal{pr}_{U}\{ \eta - s; (\mathbf{J}_{m}^\ast,s,k),\mathbf{f} \} \mu^\ast_{Y_{m}^\ast k}(s; \mathbf{J}_{m}^\ast,\boldsymbol\theta)   \diff s.
\end{align*}
We now move on to the jump-mark distribution. For $t \leq \eta$,
\begin{align*}
    P( Y_{m+1} = k \mid \mathbf{J}_{m}, T_{m+1}=t ) &= P( Y^\ast_{m+1} = k \mid \mathbf{J}_{m}^\ast, T^\ast_{m+1}=t, T^\ast_{m+1}+U^\ast_{m+1} \leq \eta ) \\
    &= \frac{ P( Y^\ast_{m+1} = k, U^\ast_{m+1} \leq \eta-t \mid \mathbf{J}_{m}^\ast, T^\ast_{m+1}=t ) }{ P( U^\ast_{m+1} \leq \eta-t  \mid \mathbf{J}_{m}^\ast, T^\ast_{m+1}=t ) }.
\end{align*}
Similar factorizations to before lead to
\begin{align*}
    P( Y_{m+1} = k \mid \mathbf{J}_{m}, T_{m+1}=t  ) &= \frac{ \textnormal{pr}_{U}\{ \eta - t; (\mathbf{J}_{m}^\ast,t,k),\mathbf{f} \} \mu^\ast_{Y_{m}^\ast k}(t; \mathbf{J}_{m}^\ast,\boldsymbol\theta)     }{ \sum_{\ell = 1}^J \textnormal{pr}_{U}\{ \eta - t; (\mathbf{J}_{m}^\ast,t,\ell),\mathbf{f} \} \mu^\ast_{Y_{m}^\ast \ell}(t; \mathbf{J}_{m}^\ast,\boldsymbol\theta)   } \\
    &= \frac{\nu_{Y_m^\ast k}(t;\mathbf{J}_m^\ast,\boldsymbol\theta, \mathbf{f})}{\nu_{Y_m^\ast \bigcdot}(t;\mathbf{J}_m^\ast,\boldsymbol\theta, \mathbf{f})}.
\end{align*} 
Therefore,
\begin{align*}
    &\frac{P(T_{m+1} \in \diff t, Y_{m+1}=k \mid \mathbf{J}_m)}{P(T_{m+1} \geq t \mid \mathbf{J}_m)} \\
    &= \frac{\exp\left\{-\int_{T_m}^t \mu_{Y_m \bigcdot}^\ast(s;\mathbf{J}_m,\boldsymbol\theta) \diff s\right\} \nu_{Y_m k}(t; \mathbf{J}_m, \boldsymbol\theta, \mathbf{f}) }{\textnormal{pr}_{U}( \eta - T_m; \mathbf{J}_{m},\mathbf{f} )-\int_{T_m}^t \exp\left\{-\int_{T_m}^s \mu_{Y_m \bigcdot}^\ast(v;\mathbf{J}_m,\boldsymbol\theta) \diff v \right\} \nu_{Y_m \bigcdot}(s;\mathbf{J}_m,\boldsymbol\theta,\mathbf{f}) \diff s } \diff t.
\end{align*}
By Jacod's formula for the intensity, see Proposition~(3.1) of~\citet{Jacod:1975}, we immediately get 
\begin{align*}
\mu_{jk}\{t; \mathbf{J}(t-), \boldsymbol\theta, \mathbf{f}\} = \gamma(t) \times
\nu_{jk}\{t; \mathbf{J}(t-), \boldsymbol\theta, \mathbf{f}\}
\end{align*}
on $\{I_j(t)=1\}$ using 
$$\exp\left(-\int_{T_m}^t \mu_{Y_m \bigcdot}^\ast(s;\mathbf{J}_m,\boldsymbol\theta) \diff s\right) = 1-\int_{T_m}^t \exp\left(-\int_{T_m}^s \mu_{Y_m \bigcdot}^\ast(v;\mathbf{J}_m,\boldsymbol\theta) \diff v\right) \mu_{Y_m \bigcdot}^\ast(s;\mathbf{J}_m,\boldsymbol\theta)  \diff s.$$ 
\end{proof}

\begin{remark}
    Expressing the thinned hazard in terms of the original hazard allows one to use observations in the thinned process to estimate the original hazard. Usually one has some expert knowledge about the form of the underlying hazard, while positing a reasonable model for the thinned hazard can be more challenging. Therefore, it becomes natural to express the thinned hazard in terms of the original hazard, as is our approach. If one had instead sought to use a data-adaptive method to estimate the thinned hazard $\mu$ directly, it would have been more natural to invert the above formula to obtain an expression for the original hazard $\mu^\ast$ as a function of $\mu$. Under sufficient smoothness conditions, one can obtain a coupled system of differential equations for the hazards $\mu^\ast_{jk}\{s;\mathbf{J}(s-),\boldsymbol\theta\}$ in terms of $\mu_{jk}\{s;\mathbf{J}(s-),\boldsymbol\theta\}$. We briefly sketch the argument. Note that, on $(T_m < \infty, Y_m=j)$,
    \begin{align*}
        &\exp\left\{-\int_{T_m}^t \mu_{j \bigcdot}(s; \mathbf{J}_m, \boldsymbol\theta, \mathbf{f}) \diff s \right\} \mu_{jk}(t; \mathbf{J}_m, \boldsymbol\theta, \mathbf{f}) \\
        &= \exp\left \{-\int_{T_{m}}^t \mu^\ast_{j \bigcdot}(s ; \mathbf{J}_{m}, \boldsymbol\theta) \diff s \right \} \mu^\ast_{jk}(t; \mathbf{J}_{m},\boldsymbol\theta) \times \frac{\textnormal{pr}_{U}\{ \eta - t; (\mathbf{J}_{m},t,k),\mathbf{f} \}}{\textnormal{pr}_U(\eta - T_m ; \mathbf{J}_{m}, \mathbf{f})}.
    \end{align*}
    This follows by casting
    $\mathbb{P}(T_{m+1} \in \diff t \mid \mathbf{J}_m) \mathbb{P}(Y_{m+1} = k \mid \mathbf{J}_m, T_{m+1}=t)$ first in terms of $\mu$, invoking the definition of $\mu_{jk}$, and secondly using the expressions obtained above. Applying the logarithm and differentiating with respect to $t$ gives a differential equation for $\mu_{jk}^\ast$. Together, these differential equations form a system of nonlinear quadratic differential equations which may be solved numerically.
\end{remark}

\begin{remark}
   If instead of imposing a parametric model, one had assumed the generalized Cox model of~\citet{Dabrowska:1997}, which allows for time-dependent covariates, the derivation of the hazard in the thinned marked point process would be unchanged, implying that $\mu_{jk}$ would no longer be on the proportional hazards form. However, $\nu_{jk}$ stays on the proportional hazards form since the baseline hazard absorbs the reporting delay distribution. When using the smoothed profile likelihood of~\citet{Dabrowska:1997} and the Poisson approximation, one may thus ignore reporting delays for inference about the parametric part of a Cox model.
\end{remark}

\subsection{Approximating the hazard}

Recall 
\begin{align*}
    \check{\gamma}(t) &= \frac{1-\int_{t-D(t-)}^t P^\ast\{s; \mathbf{J}(s-), \boldsymbol\theta\} \mu^\ast_{j\bigcdot}\{s; \mathbf{J}(s-), \boldsymbol\theta\} \, \mathrm{d}s}
{1-\int_{t-D(t-)}^t P^\ast\{s; \mathbf{J}(s-), \boldsymbol\theta\} \nu_{j\bigcdot}\{s; \mathbf{J}(s-), \boldsymbol\theta, \mathbf{f}\} \, \mathrm{d}s}.
\end{align*}
Since $\nu_{jk}\{s; \mathbf{J}(s-),\boldsymbol\theta,\mathbf{f}\} \leq \mu^\ast_{jk}\{s; \mathbf{J}(s-), \boldsymbol\theta\}$, it holds that $\check{\gamma}(t) \leq 1.$ Conversely, it holds that $\nu_{j\bigcdot}\{s; \mathbf{J}(s-),\boldsymbol\theta,\mathbf{f}\} \geq 0$ so
$\check{\gamma}(t) \geq P^\ast\{t; \mathbf{J}(t-), \boldsymbol\theta\}.$ 
These are the bounds reported in the main text. To show that $\check{\mu}_{jk}$ and $\nu_{jk}$ only differ by a second order term, impose Assumption~\ref{AssumptionSmoothHazard}. 
\begin{assumption} \label{AssumptionSmoothHazard}
    As a function of $s$, the hazards $\mu_{jk}(s;\mathbf{J},\boldsymbol\theta,\mathbf{f})$ and $\nu_{jk}(s;\mathbf{J},\boldsymbol\theta,\mathbf{f})$ are continuous and almost everywhere continuously differentiable.
\end{assumption}
Assumption~\ref{AssumptionSmoothHazard} is not used elsewhere in the paper, but similar assumptions are made when deriving asymptotic properties of the estimators. Assumption~\ref{AssumptionSmoothHazard} implies that the hazards are almost everywhere continuously differentiable so Taylor's theorem with Lagrange remainder at the point $t-D(t-)$ and some straightforward calculations yield that 
\begin{align*}
    \check{\gamma}(t) &=1+ D(t-) \times \mu^\ast_{ j \bigcdot }\{ \zeta ; \mathbf{J}( \zeta-), \boldsymbol\theta\} \times  P^\ast\{\zeta; \mathbf{J}(\zeta-), \boldsymbol\theta \} \\
    &\quad \times \bigg(  \frac{\nu_{j \bigcdot}\{\zeta;  \mathbf{J}(\zeta-), \boldsymbol\theta, \mathbf{f}\}}{\mu^\ast_{j \bigcdot}\{\zeta ;  \mathbf{J}(\zeta-), \boldsymbol\theta\}} \Big[1- \int_{t-D(t-)}^\zeta P^\ast\{s ; \mathbf{J}(s-), \boldsymbol\theta\} \mu^\ast_{j \bigcdot}\{s;  \mathbf{J}(s-), \boldsymbol\theta\} \diff s \Big]  \\
    & \qquad \qquad +\int_{t-D(t-)}^\zeta P^\ast\{s ; \mathbf{J}(s-), \boldsymbol\theta\} \nu_{j \bigcdot}\{s;  \mathbf{J}(s-), \boldsymbol\theta,\mathbf{f}\} \diff s -1 \bigg) \\ & \quad \times \bigg[1 - \int_{t-D(t-)}^\zeta P^\ast\{s ; \mathbf{J}(s-), \boldsymbol\theta\} \nu_{j \bigcdot}\{s;  \mathbf{J}(s-), \boldsymbol\theta, \mathbf{f}\} \diff s \bigg]^{-2} \\
    &= 1+R(t)
\end{align*}
for some $\zeta$ between $t-D(t-)$ and $t$. Since $\check{\gamma}$ is bounded above by one, the term $R(t)$ is negative. The term $\nu_{j \bigcdot}\{\zeta;  \mathbf{J}(\zeta-), \boldsymbol\theta, \mathbf{f}\} /\mu^\ast_{j \bigcdot}\{\zeta ;  \mathbf{J}(\zeta-), \boldsymbol\theta\} \times \Big[1- \int_{t-D(t-)}^\zeta P^\ast\{s ; \mathbf{J}(s-), \boldsymbol\theta\} \allowbreak \times \mu^\ast_{j \bigcdot}\{s;  \mathbf{J}(s-), \boldsymbol\theta\} \diff s \Big]$ is positive and $R(t)$ thus necomes more negative by removing this term. We thus get the following crude bound:
\begin{align*}
    \vert R(t) \vert &\leq 
    D(t-) \times \mu^\ast_{j \bigcdot}\{ \zeta ; \mathbf{J}(\zeta-),\boldsymbol\theta\} \times P^\ast\{\zeta; \mathbf{J}(\zeta-),\boldsymbol\theta\}  \\
    & \quad \times \frac{1-\int_{t-D(t-)}^\zeta P^\ast\{s ; \mathbf{J}(s-),\boldsymbol\theta\} \nu_{j \bigcdot}\{s;  \mathbf{J}(s-), \boldsymbol\theta,\mathbf{f}\} \diff s}{[1 - \int_{t-D(t-)}^\zeta P^\ast\{s; \mathbf{J}(s-),\boldsymbol\theta\} \nu_{j \bigcdot}\{s;  \mathbf{J}(s-), \boldsymbol\theta,\mathbf{f}\} \diff s ]^2 } \\
    &= D(t-) \times \mu^\ast_{j \bigcdot}\{ \zeta ; \mathbf{J}(\zeta-),\boldsymbol\theta\} \times \check{\gamma}(\zeta)\\
    & \leq D(t-) \times \mu^\ast_{ j \bigcdot }\{ \zeta ; \mathbf{J}(\zeta-),\boldsymbol\theta\}.
\end{align*}
The bound in the main text now follows by using that 
$$\nu_{jk}\{t;\mathbf{J}(t-),\boldsymbol\theta,\mathbf{f}\} \leq \mu^\ast_{j \bigcdot}\{t;\mathbf{J}(t-),\boldsymbol\theta\} \leq \\ \sup_{s \in [t-D(t-),t]} \mu^\ast_{j \bigcdot}\{s;\mathbf{J}(s-),\boldsymbol\theta\}.$$

\section{Fully Endogenous Reporting Delays} \label{supplement:Delays}

\subsection{Estimation procedure}

In the proposed method, adjudication enters exogenously while reporting delays enter partially exogenously and partially endogenously with respect to the marked point process of interest. Included jumps are on the form $(T_m^\ast,Y_m^\ast) \times 1_{(T_m^\ast + U_m^\ast \leq \eta)}$ instead of the fully endogenous form $(T_m^\ast,Y_m^\ast,U_m^\ast) \times 1_{(T_m^\ast + U_m^\ast \leq \eta)} $, while the adjudication information $\mathcal{A}_{\langle \eta \rangle,\eta}$ only affects the imputation of $\xi_{\langle \eta \rangle}$. Here, we investigate what happens if reporting delays instead are made fully endogenous. A result of this is that one must describe how subsampling by prior reporting delays affects the intensity of future jumps. However, if one is able to specify this correctly, for example through independence assumptions, there would be an efficiency gain in using the resulting estimator. If $\mathcal{A}_{m,t}$ is generated by a marked point process, one could similarly make the adjudication information endogenous by extending with these jump-times and jump-marks. However, it seems reasonable to retain adjudication as exogeneous, as opposed to reporting delays, since the latter often stem from the subject's conduct, while adjudication likely results from the internal processes of the data collector. This leads to a marked point process $(T^\ast_m,Y^\ast_m,U^\ast_m)_{m \geq 1}$. To relate the hazards of this model to the hazards of $(T^\ast_m,Y^\ast_m)_{m \geq 1}$, one has to posit a model for how the hazards subsampled on reporting delays compare to the ones where one does not condition on prior reporting delays. An example of such a model could be the full independence $(T^\ast_{m+1},Y^\ast_{m+1}) \indep U_{1}^\ast,\ldots,U_{m}^\ast \mid \mathbf{J}^\ast_{m}$, meaning that the hazards are unchanged by the additional information. Define $\mathbf{J}_U(t)$ as the time $t$ event history of $(T_m,Y_m,U_m)_{m \geq 1}$. Redoing the calculations of the distribution of the thinned marked point process leads to the partial likelihood
\begin{align*}
\begin{split}
\ell(\boldsymbol\theta;\mathbf{f}) &= \sum_{j,k=1}^J \int_{V}^{C} \log \mu_{jk}\{t;\mathbf{X},\mathbf{J}_U(t-),\boldsymbol\theta,\mathbf{f}\}  \diff N_{jk}(t) \\
&\quad- \int_V^{C}  I_j(t) \mu_{jk}\{t;\mathbf{X},\mathbf{J}_U(t-),\boldsymbol\theta,\mathbf{f}\} \diff t  \\
&\quad+\int_V^{C}  \log \bigg( \frac{\textnormal{pr}_{U}[\diff U_{N_{\bigcdot \bigcdot}(t)}; \{\mathbf{J}_U(t-),t,k\}, \mathbf{f}]}{\textnormal{pr}_{U}[\eta-t; \{\mathbf{J}_U(t-),t,k\},\mathbf{f}] } \bigg)  \diff N_{jk}(t). 
\end{split}
\end{align*}
In a slight abuse of notation, we have reused the previous notations $\ell(\boldsymbol\theta;\mathbf{f})$ and $\mu_{jk}$, with the latter now being $\mu_{jk}\{t; \mathbf{J}_U(t-), \boldsymbol\theta, \mathbf{f}\} = \gamma(t) \times \nu_{jk}\{t; \mathbf{J}_U(t-), \boldsymbol\theta, \mathbf{f}\}$ where $\gamma(t)$ equals
\begin{align*}
\frac{1-\int_{t-D(t-)}^t P^\ast\{s; \mathbf{J}(s-), \boldsymbol\theta\} \mu^\ast_{j\bigcdot}\{s; \mathbf{J}(s-), \boldsymbol\theta\} \, \mathrm{d}s}
{\textnormal{pr}_U[\eta - \{t-D(t-)\}; \mathbf{J}\{t-D(t-)\}, \mathbf{f}]-\int_{t-D(t-)}^t P^\ast\{s; \mathbf{J}(s-), \boldsymbol\theta\} \nu_{j\bigcdot}\{s; \mathbf{J}_U(s-), \boldsymbol\theta, \mathbf{f}\} \, \mathrm{d}s}
\end{align*}
and $\nu_{jk}\{t; \mathbf{J}_U(t-), \boldsymbol\theta, \mathbf{f}\}=\mu^\ast_{jk}\{t; \mathbf{J}(t-), \boldsymbol\theta\} \textnormal{pr}_{U}\big[\eta-t ; \{\mathbf{J}_U(t-),t,k\}, \mathbf{f}\big]$. Adjusting for adjudication then leads us to introduce $\ell(\mathbf{Z},\boldsymbol\theta;\mathbf{f},\mathbf{g}) = \mathbb{E}_\mathbf{g}\{\ell(\boldsymbol\theta;\mathbf{f}) \mid \mathbf{Z}\} = w(0,\mathbf{Z};\mathbf{g})\ell(\boldsymbol\theta;\mathbf{f},0)+w(1,\mathbf{Z};\mathbf{g})\ell(\boldsymbol\theta;\mathbf{f},1)$. The estimation algorithm becomes more involved, even under the Poisson approximation, since $\mathbf{f}$ in the optimizing step now additionally enters into the first two terms of $\ell(\boldsymbol \theta, \mathbf{f})$. For fixed $\mathbf{f}$, the estimation of $\boldsymbol\theta$ proceeds as in Section~\ref{sec:est_theta}. This suggests to optimize $\ell(\boldsymbol\theta;\mathbf{f})$ by defining $\hat{\boldsymbol\theta}_n(\mathbf{f})$ as the estimated value of $\boldsymbol\theta$ for a fixed value $\mathbf{f}$, and then optimizing $\ell\{\hat{\boldsymbol\theta}_n(\mathbf{f});\mathbf{f}\}$ over $\mathbf{f}$ using some suitable numerical optimizer. 

\subsection{Efficiency}

To discuss efficiency, we assume as in Section~\ref{sec:Asymptotics} that scores exist and we work in a neighborhood of the true parameter where the score has a unique zero. We also impose certain regularity conditions such that $\hat{\boldsymbol\theta}_n$ is asymptotically normal in both cases. In principle, the efficiency gain or loss by making the reporting delays endogenous can be inferred by comparing the expressions for the asymptotic variance of $\hat{\boldsymbol\theta}_n$ in the two cases. The relation between the covariance matrices, however, appears unclear in the general case. Since the space of distributions satisfying the moment conditions in the case where reporting delays are endogenous is strictly contained in the space of distributions satisfying the moment conditions in the partially endogenous case, the semiparametric efficiency bound in the fully endogenous case is smaller or equal to the one in the partially endogenous case. Consequently, when using the optimal weighting estimator as discussed in Section~\ref{sec:Asymptotics}, the estimator in the fully endogenous case is at least as efficient as the one in the partially endogeneous case.

\subsection{Admissibility in data application}

The assumption of full independence is not admissible for the disability semi-Markov application of Section~\ref{sec:Application}. This is due to the fact that it is assumed that death has no reporting delay, while monotone reporting remains imposed; this has the effect that knowing $(T^\ast_{m},U^\ast_{m})$ implies that a death cannot have occurred in $[T^\ast_{m},T^\ast_{m}+U^\ast_{m}]$. It might then instead be natural to say that the only extra information gained by conditioning on reporting delays is that the subject has not died before $T^\ast_{m}+U^\ast_{m}$ In Section 6 of \citet{Hoem:1972}, the author obtains an expression for the transition probabilities and rates of a semi-Markov process conditional on the time of absorption not having occurred before a given time. Let the semi-Markov transition probabilities be denoted $p^\ast_{jk}(s,t,u,z) = \mathbb{P}( I^\ast_k(t)=1, D^\ast(t) \leq z \mid I^\ast_j(s) = 1, D^\ast(s) = u)$. Assume the state space $1,\ldots,J$ can be partitioned into two parts $\mathbf{R}_1$ and $\mathbf{R}_2$, where $\mathbf{R}_2$ is absorbing in the sense that $\mathbf{R}_1$ cannot be reached from $\mathbf{R}_2$. The transition rates conditional on being in $\mathbf{R}_1$ at time $t$ are given by
\begin{align*}
\mu^\ast_{jk}(s,u;t)=\mu^\ast_{jk}(s,u)\frac{ \sum_{i \in \mathbf{R}_1} p^\ast_{ki}(s,t,0,\infty)  }{\sum_{i \in \mathbf{R}_1} p^\ast_{ji}(s,t,u,\infty) }
\end{align*}
for $s < t$ and $j \in \mathbf{R}_1$, while $\mu^\ast_{jk}(s,u;t) = \mu^\ast_{jk}(s,u)$ for $s \geq t$. So in principle, this situation may also be handled by our approach. However, it would make the estimation procedure considerably more computationally demanding since the inclusion of the transition probabilities $p^\ast_{jk}$ adds another layer of numerical integration that has to be performed separately for each candidate value of $\boldsymbol\theta$.

\section{Implementation of Estimation Procedure} \label{supplement:Implement}
A general approach to implementation of the estimation procedure for our parametric models is to discretize the integrals and use Poisson regression. For counting process likelihoods written in terms of the usual hazards, this is well-known and has been noted in, for instance,~\citet{Lindsey:1995}, but we see that this is also true for the reporting delay distribution parameterized in terms of reverse time hazards, as well as for the Poisson approximation of the likelihood. Contrary to the other supplements, here the dependence on baseline covariates is made explicit.

Select a partition $0=t_0<\ldots<t_A = \eta$ and discretize the likelihood according to this partition. Then
 \begin{align*}
    \ell_\xi(\mathbf{Z}^\textnormal{obs},\mathbf{g}) &= \sum_{m=1}^{\langle \eta \rangle} \sum_{j,k = 1}^K \sum_{a=0}^{A-1} \log \omega_{jk}(t_a;\mathcal{H}_{t_a},\mathbf{g}) O_{jk}^{\xi}(a,m) \\
    & \qquad \qquad \qquad \quad - \omega_{jk}(t_a;\mathcal{H}_{t_a},\mathbf{g}) E^\xi_{j}(a,m) + \varepsilon_\xi(A)
\end{align*}
where $O^\xi_{jk}(a,m) = N_{m,jk}(t_{a+1})-N_{m,jk}(t_{a})$ and $E^\xi_j(a,m) = \int_{t_{a}}^{t_{a+1}} I_{m,j}(s) \diff s$ are the adjudication occurrences and exposures, respectively, while $\varepsilon_\xi(A) \rightarrow 0$ for $A \rightarrow \infty$. Multiplying with the exposures inside the logarithm term, we recognize this as the likelihood corresponding to the situation where $O^\xi_{jk}(a,m)$, $(a=0,\ldots,A-1; m=1,\ldots, \langle \eta \rangle; j,k = 1,\ldots,K)$, are mutually independent Poisson random variables with mean $\log \{ \omega_{jk}(t_a;\mathcal{H}_{t_a},\mathbf{g}) E^\xi_{j}(a,m) \}$. This further corresponds to a generalized linear model with the Poisson family, log link, mean $\log  \omega_{jk}(t_a;\mathcal{H}_{t_a},\mathbf{g})$, and offset $\log E^\xi_{j}(a,m)$. The discretization may either be seen as an approximation which converges to the true criterion function for $A \rightarrow \infty$, or it may be seen as the true likelihood if the partition is chosen to be equal to the precision of the time measurement. 

Similarly for the reporting delays, let $\overline{t}_a=(t_{a+1}+t_a)/2$ and note that
\begin{align*}
    \ell_U(\mathbf{f}) &= \sum_{m=1}^{N^c_{\bigcdot \bigcdot}(\eta)} \sum_{a=0}^{A-1} \log \alpha \{ \overline{t}_a;\mathbf{J}^c(T^c_m),\mathbf{X},\mathbf{f} \} O^U(a,m) \\
    & \qquad \qquad \quad  - \alpha \{ \overline{t}_a;\mathbf{J}^c(T^c_m),\mathbf{X},\mathbf{f} \} E^U(a,m) + \varepsilon_U(A),
\end{align*}
where $O^U(a,m) = 1_{(U^c_m \in (t_{a},t_{a+1}] )}$ and $E^U(a,m) = \int_{t_{a}}^{t_{a+1}} 1_{(U^c_m \leq s \leq \eta-T^c_m)} \diff s$ are the reporting delay occurrences and exposures, while $\varepsilon_U(A) \rightarrow 0$ for $A \rightarrow \infty$. Note that the partition here is allowed to differ from the one of $\ell_\xi(\mathbf{Z}^\textnormal{obs},\mathbf{g})$ despite using the same notation. Here, the midpoints are chosen to take into account that the reverse time hazard takes the value $+\infty$ at time $0$.

For the approximate likelihood of $\boldsymbol\theta$, we get
\begin{align*}
    \ell^{\textnormal{app}}(\boldsymbol\theta;\mathbf{f}) &=  \sum_{j,k=1}^J \sum_{a=0}^{A-1} \log \mu^\ast_{jk}\{t_a ; \mathbf{J}( t_a),\mathbf{X},\boldsymbol\theta\}  O_{jk}(a) \\
& \qquad \qquad \quad - \mu^\ast_{jk}\{t_a ; \mathbf{J}(t_a),\mathbf{X}, \boldsymbol\theta\} \mathbb{P}_{U}[\eta-t_a ; \{\mathbf{J}(t_a),t_a,k\},\mathbf{X}, \mathbf{f}]  E_j(a) + \varepsilon(A)
\end{align*}
where $O_{jk}(a) = N^c_{jk}(t_{a+1})-N^c_{jk}(t_{a})$ and $E_j(a) = \int_{t_{a}}^{t_{a+1}} 1_{(V,C]}(s) I_j(s) \diff s$ are the occurrences and exposures, respectively, while $\varepsilon(A) \rightarrow 0$ for $A \rightarrow \infty$. When computing $\ell^\textnormal{app}(\mathbf{Z},\boldsymbol\theta;\mathbf{f},\mathbf{g})$, the occurrence and exposure contributions have to be weighted with the adjudication probabilities, but otherwise nothing changes. Most generalized linear model software packages allow the user to specify such weights.  

The link to Poisson regression simplifies practical implementation significantly. Another useful consequence is that the link provides a natural way to construct plots in which the fitted rates can be compared with empirical rates. If the hazards are assumed to be piecewise constant, the maximum likelihood estimates become the well-known occurrence-exposure rates. It is these occurrence-exposure rates which we refer to as the empirical rates. The piecewise constant hazard assumption can be imposed for any continuous covariate, for example calendar time, age, or the duration since the occurrence of given event. By comparing the empirical rates with the predicted rates from the fitted parametric model, one obtains a goodness-of-fit plot.

In our numerical study, that data is on a jump-time and jump-mark format. There we forego the Poisson regression, opting instead to utilize the specific choice of parametric model to obtain a direct implementation; this includes the use of a standard optimizer to maximize the resulting expressions. In our data application, the data is on a monthly grid format, and there is considerably more data, which leads us to use the Poisson regression approach outlined above, since estimation might then exactly be based on aggregated data. When the number of observations is significantly larger than the covariate space, which is the case in our data application, this saves considerable amounts of storage and computation time.

\section{Asymptotics for the Parameters} \label{supplement:Asymptotics}

Consistency of the nuisance parameter estimators $\hat{\mathbf{f}}_n$ and $\hat{\mathbf{g}}_n$ is needed in order to achieve consistency of $\hat{\boldsymbol\theta}_n$ and these estimators generally need a convergence rate of a higher order than $n^{1/4}$ in order to obtain asymptotic normality of $\hat{\boldsymbol\theta}_n$, see e.g.\ Theorem 2 of~\citet{Chen.etal:2003}. In this framework, the convergence rate is obtained with plenty to spare since parametric models under regularity conditions and correct specifications admit $n^{1/2}$-rates.

\subsection{Consistency}

Before showing asymptotic normality, we establish consistency. We focus on weak consistency, as this suffices for asymptotic distribution theory. Showing strong consistency of $\hat{\mathbf{g}}_n$, $\hat{\mathbf{f}}_n$, and $\hat{\boldsymbol\theta}_n$ would, however, require almost no additional effort. For Theorem 2.1 in~\citet{Newey:McFadden:1994}, one would change~(iv) to hold almost surely in order to obtain strong consistency. Similarly, for Theorem 1 in~\citet{Delsol:Van:2020}, one could change~(A1) such that the estimator is required to be the actual maximizer of the objective function, which is satisfied for our estimators, and further impose that~(A3) and~(A4) hold almost surely instead of in probability. Since we show~(iv) and~(A4) by a Glivenko--Cantelli class argument, which implies almost sure convergence, we could equally well claim strong consistency. 

Our approach to showing consistency and asymptotic normality is not to find minimal assumptions, but rather to formulate sufficient assumptions that are easy to check in a given application and, in particular, are not too restrictive for our applications. 

\begin{assumption} \label{assumption:g} \leavevmode
    \begin{enumerate}[(i)]
        \item Different $\mathbf{g}$ give rise to $\omega_{jk}$ that are not almost surely equal.
        \item The $\mathbf{g}$ parameter space $\mathbb{G}$ is a compact subset of a Euclidean space, and $\mathbf{g}_0$ belongs to the interior.
        \item The hazards $\omega_{jk}(t;\mathbf{H},\mathbf{g})$ are continuous in $\mathbf{g}$.
        \item $\sup_{\mathbf{g}} \vert \ell_\xi(\mathbf{Z}^\textnormal{obs},\mathbf{g}) \vert$ has finite expectation and $\sup_{\mathbf{g}, \mathbf{H}}\omega_{jk}(t;\mathbf{H},\mathbf{g})$ are Lebesgue-integrable on $[0,\eta]$.
    \end{enumerate}
\end{assumption}

In the following, we use the usual stochastic order notation $o_P$, that is we write $x_n = o_P(r_n)$ if $\lvert x_n \rvert /r_n \rightarrow 0$ in probability. 

\begin{proposition} \label{proposition:gConsistency}
    Under Assumptions~\ref{AssumptionMonotoneReporting}-\ref{assumption:oneAdj} and Assumption~\ref{assumption:g}, it holds that $\norm{\hat{\mathbf{g}}_n- \mathbf{g}_0} = o_P(1)$.
\end{proposition}
\begin{proof}
We verify the conditions of Theorem 2.1 in~\citet{Newey:McFadden:1994}. Condition~(i) is satisfied by non-negativity of the Kullback-Leibler divergence, also known as Gibbs' inequality or the information inequality, and the identifiability of the model implied by Assumption~\ref{assumption:g}(i); each contribution to the sum in $\ell_\xi$ is a usual log-likelihood of a marked point process which is maximized in $\mathbf{g}_0$ by Gibbs' inequality, and the identifiability ensures that a different $\mathbf{g}$ would lead to a smaller value for at least one of the contributions using Gibbs' inequality and the one-to-one correspondence between compensators and distributions. Condition~(ii) is satisfied by Assumption~\ref{assumption:g}(ii). We now consider Condition~(iii).  For a given 
$\mathbf{g} \in \mathbb{G}$, take $(\mathbf{g}_j)_{j \geq 1}$ with $\mathbf{g}_j \rightarrow \mathbf{g}$. By Assumption~\ref{assumption:g}(iv), we may use dominated convergence which implies $\lim_{j \rightarrow \infty}\mathbb{E}\{\ell_{\xi}(\mathbf{Z}^{\textnormal{obs}},\mathbf{g}_j)\}=\mathbb{E}\{\lim_{j \rightarrow \infty}\ell_{\xi}(\mathbf{Z}^{\textnormal{obs}},\mathbf{g}_j)\}$. Again by Assumption~\ref{assumption:g}(iv), dominated convergence and Assumption~\ref{assumption:g}(iii) gives $\lim_{j \rightarrow \infty}\ell_{\xi}(\mathbf{Z}^{\textnormal{obs}},\mathbf{g}_j) = \ell_{\xi}(\mathbf{Z}^{\textnormal{obs}},\mathbf{g})$. Thus, Condition~(iii) holds. Condition~(iv) corresponds to showing that $\{ \ell_{\xi}(\cdot,\mathbf{g}) : \mathbf{g} \in \mathbb{G} \}$ is Glivenko--Cantelli. We will use Example~19.8 of~\citet{Van:1998}. We note that $\mathbf{g} \mapsto \ell_{\xi}(\mathbf{Z}^\textnormal{obs},\mathbf{g})$ is continuous for any $\mathbf{Z}^\textnormal{obs}$ by previous observations and has an integrable envelope function by Assumption~\ref{assumption:g}(iv). Hence the class is Glivenko--Cantelli and Condition~(iv) is satisfied.
\end{proof}

\begin{assumption} \label{assumption:f} \leavevmode
    \begin{enumerate}[(i)]
        \item Different $\mathbf{f}$ gives rise to $\alpha$ that are not almost surely equal.
        \item The $\mathbf{f}$ parameter set $\mathbb{F}$ is a compact subset of a Euclidean space, and $\mathbf{f}_0$ belongs to the interior.
        \item The reverse time hazard $\alpha(t;\mathbf{J},\mathbf{X},\mathbf{f})$ is continuous in $\mathbf{f}$.
        \item $\sup_{\mathbf{f}} \vert \ell_U(\mathbf{f};k) \vert$ for $k=0,1$ has finite expectation and $\sup_{\mathbf{f},\mathbf{J},\mathbf{X}}\alpha(t;\mathbf{J},\mathbf{X},\mathbf{f})$ is Lebesgue-integrable on $[0,\eta]$.
        \item The adjudication probabilities $w(1,\mathbf{Z};\mathbf{g})$ are continuous in $\mathbf{g}$.
    \end{enumerate}
\end{assumption}

\begin{proposition} \label{proposition:fConsistency}
    Under Assumptions~\ref{AssumptionMonotoneReporting}-\ref{assumption:oneAdj}, Assumption~\ref{assumption:f}, and if $\norm{\hat{\mathbf{g}}_n-\mathbf{g}_0} = o_P(1)$, it holds that $\Vert \hat{\mathbf{f}}_n - \mathbf{f}_0 \Vert = o_P(1)$.
\end{proposition}
\begin{proof}
We verify the conditions in Theorem 1 of~\citet{Delsol:Van:2020}. Condition~(A1) is satisfied by the specification of $\hat{\mathbf{f}}_n$ as the maximizer of the objective function. By standard arguments, Condition~(A2) is satisfied if $\mathbb{E}\{\ell_U(\mathbf{Z},\mathbf{f};\mathbf{g}_0)\}$ is continuous in $\mathbf{f}$ since $\mathbb{F}$ is compact. This is because $\mathbb{F}_\varepsilon = \{\mathbf{f} \in \mathbb{F} : \norm{\mathbf{f}-\mathbf{f}_0} \geq \varepsilon\}$ is closed and hence compact, so in this case there exists $\mathbf{f}_\varepsilon$ satisfying $\sup_{\mathbf{f} \in \mathbb{F}_\varepsilon} \mathbb{E}\{\ell_U(\mathbf{Z},\mathbf{f};\mathbf{g}_0)\} = \mathbb{E}\{\ell_U(\mathbf{Z},\mathbf{f}_\varepsilon;\mathbf{g}_0)\} < \mathbb{E}\{\ell_U(\mathbf{Z},\mathbf{f}_0;\mathbf{g}_0)\}$ by uniqueness of the maximum implied by Gibbs' inequality and the identifiability implied by Assumption~\ref{assumption:f}(i). We show simultaneous continuity in $\mathbf{f}$ and $\mathbf{g}$, as this is useful for checking other regularity conditions. The continuity is a direct consequence of Assumptions~\ref{assumption:f}(iii)-(v) by two applications of dominated convergence since $\sup_{\mathbf{f}, \mathbf{g}} \vert \ell_U(\mathbf{Z},\mathbf{f};\mathbf{g})\vert \leq \sup_{\mathbf{f}} \vert \ell_U(\mathbf{f};1)\vert + \sup_{\mathbf{f}} \vert \ell_U(\mathbf{f};0)\vert$. Condition~(A3) follows by the assumption that $\hat{\mathbf{g}}_n$ is weakly consistent. Condition~(A4) follows if the class of functions $\{ \ell_U(\cdot,\mathbf{f};\mathbf{g}) : \mathbf{f} \in \mathbb{F}, \mathbf{g} \in \mathbb{G}\}$ is Glivenko--Cantelli. We use Example 19.8 in~\citet{Van:1998}. Continuity and the existence of an integrable envelope function for the class follows by previous observations. Hence Condition~(A4) is satisfied. For Condition~(A5), we take $\mathbf{g}_j \rightarrow \mathbf{g}_0$ for $j \rightarrow \infty$. Note that by the triangle inequality, $\vert \mathbb{E}\{\ell_U(\mathbf{Z},\mathbf{f};\mathbf{g}_j)\}-\mathbb{E}\{\ell_U(\mathbf{Z},\mathbf{f};\mathbf{g}_0)\} \vert \leq \mathbb{E}\{\vert w( 0,\mathbf{Z}; \mathbf{g}_j) - w( 0,\mathbf{Z}; \mathbf{g}_0) \vert \times \sup_{\mathbf{f}}\vert\ell_U(\mathbf{f};0)\vert + \vert w( 1,\mathbf{Z}; \mathbf{g}_j) - w( 1,\mathbf{Z}; \mathbf{g}_0) \vert \times \sup_{\mathbf{f}} \vert \ell_U(\mathbf{f};1)\vert\}$. The right-hand side does not depend on $\mathbf{f}$ and converges to $0$ for $\mathbf{g}_j \rightarrow \mathbf{g}_0$ by dominated convergence and Assumption~\ref{assumption:f}(v), hence Condition~(A5) holds.
\end{proof}

\begin{assumption} \label{assumption:theta} \leavevmode
    \begin{enumerate}[(i)]
        \item Different $\boldsymbol\theta$ gives rise to $\mu^\ast_{jk}$ that are not almost surely equal.
        \item The $\boldsymbol\theta$ parameter space $\Theta$ is compact, and $\boldsymbol\theta_0$ belongs to the interior.
        \item $\sup_{\boldsymbol\theta, \mathbf{f}, \mathbf{g}} \vert \ell(\mathbf{Z},\boldsymbol\theta;\mathbf{f},\mathbf{g}) \vert$ has finite expectation.
        \item $\sup_{\boldsymbol\theta, \mathbf{J}, \mathbf{X}}\mu^\ast_{jk}(t;\mathbf{J},\mathbf{X},\boldsymbol\theta)$ is Lebesgue-integrable on $[0,\eta]$.
        \item The hazards $\mu^\ast_{jk}(t;\mathbf{J},\mathbf{X},\boldsymbol\theta)$ are continuous in $\boldsymbol\theta$.
        \item The partial derivatives of $\ell(\mathbf{Z},\boldsymbol\theta;\mathbf{f},\mathbf{g})$ with respect to the coordinates of $\boldsymbol\theta$, $\mathbf{f}$, and $\mathbf{g}$ all exist and are continuous. Furthermore, $\sup_{\boldsymbol\theta,\mathbf{f},\mathbf{g}} \Vert \nabla_{(\boldsymbol\theta,\mathbf{f},\mathbf{g})} \ell(\mathbf{Z},\boldsymbol\theta;\mathbf{f},\mathbf{g}) \Vert$ has finite expectation.
    \end{enumerate}
\end{assumption}

\begin{proposition} \label{proposition:thetaConsistent}
    Under Assumptions~\ref{AssumptionMonotoneReporting}-\ref{assumption:oneAdj}, Assumption~\ref{assumption:theta}, and if $\Vert \hat{\mathbf{g}}_n - \mathbf{g}_0 \Vert = o_P(1)$ and $\Vert \hat{\mathbf{f}}_n - \mathbf{f}_0 \Vert = o_P(1)$, it holds that $\Vert \hat{\boldsymbol\theta}_n -\boldsymbol\theta_0 \Vert=o_P(1)$.
\end{proposition}
\begin{proof}
    We verify the conditions in Theorem 1 of~\citet{Delsol:Van:2020}. Condition~(A1) is satisfied by the specification of $\hat{\theta}_n$ as the maximizer of the objective function. Similarly to Proposition~\ref{proposition:fConsistency}, Condition~(A2) is satisfied by standard arguments if $E\{\ell(\mathbf{Z},\boldsymbol\theta;\mathbf{g}_0,\mathbf{f}_0)\}$ is continuous in $\boldsymbol\theta$ since $\Theta$ is compact by Assumption~\ref{assumption:theta}(ii) and the model is identifiable by Assumption~\ref{assumption:theta}(i). We start by showing the Lipschitz-type property $\vert \ell(\mathbf{Z},\boldsymbol\theta_1;\mathbf{f}_1,\mathbf{g}_1) - \ell(\mathbf{Z},\boldsymbol\theta_2;\mathbf{f}_2,\mathbf{g}_2) \vert \leq b(\mathbf{Z}) \norm{(\boldsymbol\theta_1,\mathbf{f}_1,\mathbf{g}_1)-(\boldsymbol\theta_2,\mathbf{f}_2,\mathbf{g}_2)}$ with $b$ being an integrable function, as this is useful for showing other regularity conditions. If this holds, we see  that
    \begin{align*}
        \Big\vert \mathbb{E}\{\ell(\mathbf{Z},\boldsymbol\theta_1;\mathbf{f}_1,\mathbf{g}_1)\} - \mathbb{E}\{\ell(\mathbf{Z},\boldsymbol\theta_2;\mathbf{f}_2,\mathbf{g}_2)\} \Big\vert &\leq  \mathbb{E}\{\vert \ell(\mathbf{Z},\boldsymbol\theta_1;\mathbf{f}_1,\mathbf{g}_1) 
 - \ell(\mathbf{Z},\boldsymbol\theta_2;\mathbf{f}_2,\mathbf{g}_2)\vert\} \\
 &\leq \mathbb{E}\{b(\mathbf{Z}) \norm{(\boldsymbol\theta_1,\mathbf{f}_1,\mathbf{g}_1)-(\boldsymbol\theta_2,\mathbf{f}_2,\mathbf{g}_2)}\}
    \end{align*}
and since $b$ is integrable, dominated convergence implies that the right-hand side converges to $0$ for $(\boldsymbol\theta_1,\mathbf{f}_1,\mathbf{g}_1) \rightarrow (\boldsymbol\theta_2,\mathbf{f}_2,\mathbf{g}_2)$, and hence $\mathbb{E}\{\ell(\mathbf{Z},\boldsymbol\theta;\mathbf{f},\mathbf{g})\}$ is continuous in $\boldsymbol\theta$, $\mathbf{f}$, and $\mathbf{g}$. To show the Lipschitz-type property, note that Assumption~\ref{assumption:theta}(vi) implies that $\ell(\mathbf{Z},\boldsymbol\theta;\mathbf{f},\mathbf{g})$ is differentiable in $(\boldsymbol\theta,\mathbf{f},\mathbf{g})$ for any fixed $\mathbf{Z}$. The mean-value theorem thus implies 
\begin{align*}
&\ell(\mathbf{Z},\boldsymbol\theta_1;\mathbf{f}_1,\mathbf{g}_1) 
 - \ell(\mathbf{Z},\boldsymbol\theta_2;\mathbf{f}_2,\mathbf{g}_2) \\
 &= \nabla_{(\boldsymbol\theta,\mathbf{f},\mathbf{g})} \ell[\mathbf{Z}, \{1-k(\mathbf{Z})\}(\boldsymbol\theta_1;\mathbf{f}_1,\mathbf{g}_1)+k(\mathbf{Z})(\boldsymbol\theta_2;\mathbf{f}_2,\mathbf{g}_2) ]^\top \{(\boldsymbol\theta_1,\mathbf{f}_1,\mathbf{g}_1)-(\boldsymbol\theta_2,\mathbf{f}_2,\mathbf{g}_2)\}    
\end{align*}
for some $k(\mathbf{Z}) \in (0,1)$. The Cauchy-Schwarz inequality now gives the Lipschitz-type property, with the integrability condition holding due to Assumption~\ref{assumption:theta}(vi). Hence Condition~(A2) holds. Condition~(A3) follows by the assumption that $\hat{\mathbf{g}}_n$ and $\hat{\mathbf{f}}_n$ are weakly consistent. Condition~(A4) follows if $\{ \ell(\cdot,\boldsymbol\theta ;\mathbf{f},\mathbf{g}) : \mathbf{f} \in \mathbb{F}, \mathbf{g} \in \mathbb{G}, \boldsymbol\theta \in \Theta \}$ is Glivenko--Cantelli. This is the case by Example~19.7 and Theorem~19.4 in~\citet{Van:1998} due to the Lipchitz-type property. Hence Condition~(A4) follows. Condition~(A5) also follows from the Lipchitz-type property and Assumption~\ref{assumption:theta}(ii). 
\end{proof}

We now also consider consistency of Efron's `simple' nonparametric bootstrap introduced in~\citet{Efron:1979} and described in Section~\ref{sec:Asymptotics}. Based on the previous results, we immediately get weak consistency of the bootstrap estimator. The proposition is comparable with Proposition 1 in~\citet{Hahn:1996}.
\begin{proposition} \label{proposition:bootstrapConsistency}
 Under Assumptions~\ref{AssumptionMonotoneReporting}-\ref{assumption:oneAdj} and Assumptions~\ref{assumption:g}-\ref{assumption:theta}, one has $\norm{\hat{\mathbf{g}}^\textnormal{boot}_n -\mathbf{g}_0}=o_P(1)$, $\Vert \hat{\mathbf{f}}^\textnormal{boot}_n -\mathbf{f}_0 \Vert=o_P(1)$, and $\Vert \hat{\boldsymbol\theta}^\textnormal{boot}_n -\boldsymbol\theta_0 \Vert=o_P(1)$.
\end{proposition}
\begin{proof}
     By the Glivenko--Cantelli property shown in Propositions~\ref{proposition:gConsistency}-\ref{proposition:thetaConsistent}, each of the function classes $\{ \ell_{\xi}(\cdot,\mathbf{g}) : \mathbf{g} \in \mathbb{G} \}$, $\{ \ell_U(\cdot,\mathbf{f};\mathbf{g}) : \mathbf{f} \in \mathbb{F}, \mathbf{g} \in \mathbb{G}\}$, and $\{ \ell(\cdot,\boldsymbol\theta;\mathbf{f},\mathbf{g}) : \mathbf{f} \in \mathbb{F}, \mathbf{g} \in \mathbb{G}, \boldsymbol\theta \in \Theta \}$ satisfy the uniform convergence $\sup_{q} \norm{n^{-1}\sum_{i=1}^n q(\mathbf{Z}^{\textnormal{obs}}_{i})- E[q(\mathbf{Z}^{\textnormal{obs}})] } = o_P(1)$, where $q$ runs over the function class in question. Furthermore, we have shown that there exists an integrable envelope function for each of the classes. The bootstrap Glivenko--Cantelli theorem, Theorem 2.6 in~\citet{Gine:Zinn:1990}, then implies that 
     $$\sup_{q} \norm{n^{-1}\sum_{i=1}^n q(\mathbf{Z}^{\textnormal{obs}}_{i})- n^{-1}\sum_{i=1}^n q(\mathbf{Z}^{\textnormal{obs}}_{ni}) } = o_P(1)$$ 
     almost surely. Consequently, the bootstrap estimators maximize their respective non-bootstrap empirical objective functions up to an $o_P(1)$-term, which is sufficient to invoke both Theorem 2.1 of~\citet{Newey:McFadden:1994} and Theorem~1 of~\citet{Delsol:Van:2020}.
\end{proof}
A different approach to proving Proposition~\ref{proposition:bootstrapConsistency} would be to verify the conditions of Theorem~3.5 in~\citet{Arcones:Gine:1992}. 

\subsection{Asymptotic Normality}

We now consider the stacked estimating equation
\begin{align*}
n^{-1} \sum_{i=1}^{n} \mathbf{h}(\mathbf{Z}_i^{\textnormal{obs}},\boldsymbol\theta,\mathbf{f},\mathbf{g}) = \mathbf{0}
\end{align*}
for $\mathbf{h}(\mathbf{Z}^{\textnormal{obs}},\boldsymbol\theta,\mathbf{f},\mathbf{g})=\{\frac{\diff}{\diff \mathbf{g}}\ell_\xi(\mathbf{Z}^{\textnormal{obs}},\mathbf{g}), \frac{\diff}{\diff \mathbf{f}} \ell_U(\mathbf{Z},\mathbf{f};\mathbf{g}), \frac{\diff}{\diff \boldsymbol\theta}\ell(\mathbf{Z},\boldsymbol\theta;\mathbf{f},\mathbf{g})\}^\top$. Define also the part pertaining to the nuisance parameters $\mathbf{k}(\mathbf{Z}^{\textnormal{obs}},\mathbf{f},\mathbf{g})=\{\frac{\diff}{\diff \mathbf{g}}\ell_\xi(\mathbf{Z}^{\textnormal{obs}},\mathbf{g}), \frac{\diff}{\diff \mathbf{f}} \ell_U(\mathbf{Z},\mathbf{f};\mathbf{g})\}^\top$. The existence of the scores is ensured by Assumption~\ref{assumption:theta}(vi) and Assumption~\ref{assumption:asympNorm}(i).  
Introduce also
\begin{align*}
    \mathbf{H}_1 &= \mathbb{E} \{ \nabla_{\boldsymbol\theta} \mathbf{h}(\mathbf{Z}^{\textnormal{obs}},\boldsymbol\theta_0,\mathbf{f}_0,\mathbf{g}_0) \}, &  
    &\mathbf{H}_2 = \mathbb{E} \{ \nabla_{(\mathbf{f},\mathbf{g})} \mathbf{h}(\mathbf{Z}^{\textnormal{obs}},\boldsymbol\theta_0,\mathbf{f}_0,\mathbf{g}_0) \}, \\
    \mathbf{K} &= \mathbb{E} \{\nabla_{(\mathbf{f},\mathbf{g})} \mathbf{k}(\mathbf{Z}^{\textnormal{obs}},\mathbf{f}_0,\mathbf{g}_0)\},& &\boldsymbol\psi(\mathbf{Z}^{\textnormal{obs}}) = - \mathbf{K}^{-1} \mathbf{k}(\mathbf{Z}^{\textnormal{obs}},\mathbf{f}_0,\mathbf{g}_0).
\end{align*} 

\begin{assumption} \label{assumption:asympNorm} \leavevmode
    \begin{enumerate}[(i)]
        \item There exists a neighborhood of the true parameter denoted $(\Theta \times \mathbb{F} \times \mathbb{G})_\delta = \{ (\boldsymbol\theta,\mathbf{f},\mathbf{g}) : \norm{(\boldsymbol\theta,\mathbf{f},\mathbf{g})-(\boldsymbol\theta_0,\mathbf{f}_0,\mathbf{g}_0)} \leq \delta\}$ such that $\frac{\diff}{\diff \mathbf{g}}\ell_\xi(\mathbf{Z}^{\textnormal{obs}},\mathbf{g})$ and $\frac{\diff}{\diff \mathbf{f}} \ell_U(\mathbf{Z},\mathbf{f};\mathbf{g})$ exist and such that $\mathbb{E}\{ \mathbf{h}(\mathbf{Z}^{\textnormal{obs}},\boldsymbol\theta,\mathbf{f},\mathbf{g})\} = \mathbf{0}$ has a unique solution in this neighborhood.
        \item The score $\mathbf{h}(\mathbf{Z}^{\textnormal{obs}},\boldsymbol\theta,\mathbf{f},\mathbf{g})$ is dominated by a square-integrable function of $\mathbf{Z}^{\textnormal{obs}}$.
        \item The gradient of the score $\nabla_{(\boldsymbol\theta,\mathbf{f},\mathbf{g})}\mathbf{h}(\mathbf{Z}^{\textnormal{obs}},\boldsymbol\theta,\mathbf{f},\mathbf{g})$ exists and is continuous in the parameters. Furthermore, it is dominated by a square-integrable function of $\mathbf{Z}^{\textnormal{obs}}$.
        \item $(\mathbf{H}_1,\mathbf{H}_2)^T(\mathbf{H}_1,\mathbf{H}_2)$ is nonsingular.
    \end{enumerate}
\end{assumption}

We may now show asymptotic normality. For the ordinary estimator, one could proceed via Theorem~6.1 in~\citet{Newey:McFadden:1994}, but we instead use Theorem~1 of~\citet{Hahn:1996} as this also gives asymptotic normality of the bootstrap estimator. Similar conditions to those of~\citet{Hahn:1996} are found in~\citet{Arcones:Gine:1992}.

\begin{proposition} \label{prop:asympNorm}
    Under Assumptions~\ref{AssumptionMonotoneReporting}-\ref{assumption:oneAdj} and Assumptions~\ref{assumption:g}-\ref{assumption:asympNorm}, it holds that 
    \begin{enumerate}[(i)]
        \item The ordinary estimators $\hat{\mathbf{g}}_n$, $\hat{\mathbf{f}}_n$, and $\hat{\boldsymbol\theta}_n$ are asymptotically normal.
        \item $\sqrt{n}(\hat{\boldsymbol\theta}_n-\boldsymbol\theta_0) \rightarrow N(0,\mathbf{V})$ in distribution for $$V=\mathbf{H}_1^{-1} \mathbb{E}[\{\mathbf{h}(\mathbf{Z}^{\textnormal{obs}},\boldsymbol\theta_0,\mathbf{f}_0,\mathbf{g}_0)+\mathbf{H}_2 \boldsymbol\psi(\mathbf{Z}^{\textnormal{obs}})\}\{\mathbf{h}(\mathbf{Z}^{\textnormal{obs}},\boldsymbol\theta_0,\mathbf{f}_0,\mathbf{g}_0)+\mathbf{H}_2 \boldsymbol\psi(\mathbf{Z}^{\textnormal{obs}})\}^\top] (\mathbf{H}_1^{-1})^\top.$$
        \item The bootstrap estimators $\hat{\mathbf{g}}^\textnormal{boot}_n$,  $\hat{\mathbf{f}}^\textnormal{boot}_n$, and  $\hat{\boldsymbol\theta}^\textnormal{boot}_n$ are asymptotically normal and, in particular, $\sqrt{n}(\hat{\boldsymbol\theta}^\textnormal{boot}_n-\hat{\boldsymbol\theta}_n) \rightarrow N(0,\mathbf{V})$ in distribution.
    \end{enumerate}  
\end{proposition}
\begin{proof}
Due to Assumption~\ref{assumption:asympNorm}(i), our estimator can be formulated as a generalized method of moment estimator with identity weighting matrix and parameter space $(\Theta \times \mathbb{F} \times \mathbb{G})_\delta$. Theorem 1 of~\citet{Hahn:1996} therefore becomes applicable. Condition~(i) is satisfied since a) the score equation has the true parameters as its unique solution by Assumption~\ref{assumption:asympNorm}(i), b) our data $(\mathbf{Z}^\textnormal{obs}_i)_{i=1}^n$ is assumed to be independent and identically distributed, c) the solution is well-separated by standard arguments since $(\Theta \times \mathbb{F} \times \mathbb{G})_\delta$ is compact and $\mathbb{E}\{\mathbf{h}(\mathbf{Z}^{\textnormal{obs}},\boldsymbol\theta,\mathbf{f},\mathbf{g})\}$ is continuous in the parameters by dominated convergence and Assumptions~\ref{assumption:asympNorm}(ii)-(iii), d) Glivenko--Cantelli holds for the class $\{ \mathbf{h}(\cdot,\boldsymbol\theta;\mathbf{f},\mathbf{g}) : (\boldsymbol\theta,\mathbf{f},\mathbf{g}) \in (\Theta \times \mathbb{F} \times \mathbb{G})_\delta \}$ by Example~19.8 of~\citet{Van:1998} since $\mathbf{h}(\mathbf{Z}^{\textnormal{obs}},\boldsymbol\theta,\mathbf{f},\mathbf{g})$ is continuous in the parameters by Assumption~\ref{assumption:asympNorm}(iii) and has an integrable envelope function by Assumption~\ref{assumption:asympNorm}(ii), e) there is an integrable envelope function for the class as noted before, and f) the weighting matrix is equal to the identity weighting matrix. Condition~(ii) holds by specification of the estimators. Condition~(iii) holds by Assumptions~\ref{assumption:asympNorm}(ii)-(iii) via dominated convergence, since these assumptions imply that $\mathbf{h}$ is continuous in the parameters and that the integrand is dominated by an integrable function, just use $(a - b)^2 \leq 2(a^2 + b^2)$. Note that an exponent of 2 seems to be missing in the~\citet{Hahn:1996} paper; for comparison, see Theorem~21 (functional central limit theorem) in Chapter VII of~\citet{Pollard:1984}. Condition~(iv) holds if the class $\{ \mathbf{h}(\cdot,\boldsymbol\theta,\mathbf{f},\mathbf{g}) : (\boldsymbol\theta,\mathbf{f},\mathbf{g}) \in (\Theta \times \mathbb{F} \times \mathbb{G})_\delta \}$ is Donsker since Condition~(iii) holds at any parameter value, not just the true value. This can be shown via Example~19.7 and Theorem 19.5 in~\citet{Van:1998} in the same way as how the Lipschitz-type property was shown in the proof of Proposition~\ref{proposition:thetaConsistent}: It is sufficient that $\mathbf{h}$ is differentiable in the parameters and that the norm of the gradient is dominated by a square-integrable function, and these properties are ensured by Assumption~\ref{assumption:asympNorm}(iii). Condition~(v) holds by Assumptions~\ref{assumption:asympNorm}(iii)-(iv). The first part of Condition~(vi) holds by Assumption~\ref{assumption:asympNorm}(ii) and the second part holds by Theorem~3.1 of~\citet{Andersen:Dobric:1987} and their following remark when setting $p=2$ and with their function $f$ being the score $\mathbf{h}$. This is due to the following observations. Their equation~(3.1.2) holds by Assumption~\ref{assumption:asympNorm}(ii) and the triangle inequality, their equation (3.1.3) holds since the scores are both bounded and uniformly continuous with probability one by Assumptions~\ref{assumption:asympNorm}(ii)-(iii) combined with the Heine-Cantor theorem, and the parameter set is compact (and hence totally bounded) with respect to the original Euclidean metric. 
     
The variance matrix calculation is now identical to that of Theorem~6.1 in~\citet{Newey:McFadden:1994}, which leads to the desired result.
\end{proof}
Theorem~\ref{thm:asymp} is now a consequence of Propositions~\ref{proposition:gConsistency}-\ref{prop:asympNorm}. Note that the proof  can be adapted to the case where one uses the optimal weight matrix. 

\section{Additional Details for Numerical Study} \label{supplement:Numerical}

\subsection{Data-generating Process}

\begin{figure}[!ht]
    \centering
    \includegraphics[scale=0.65]{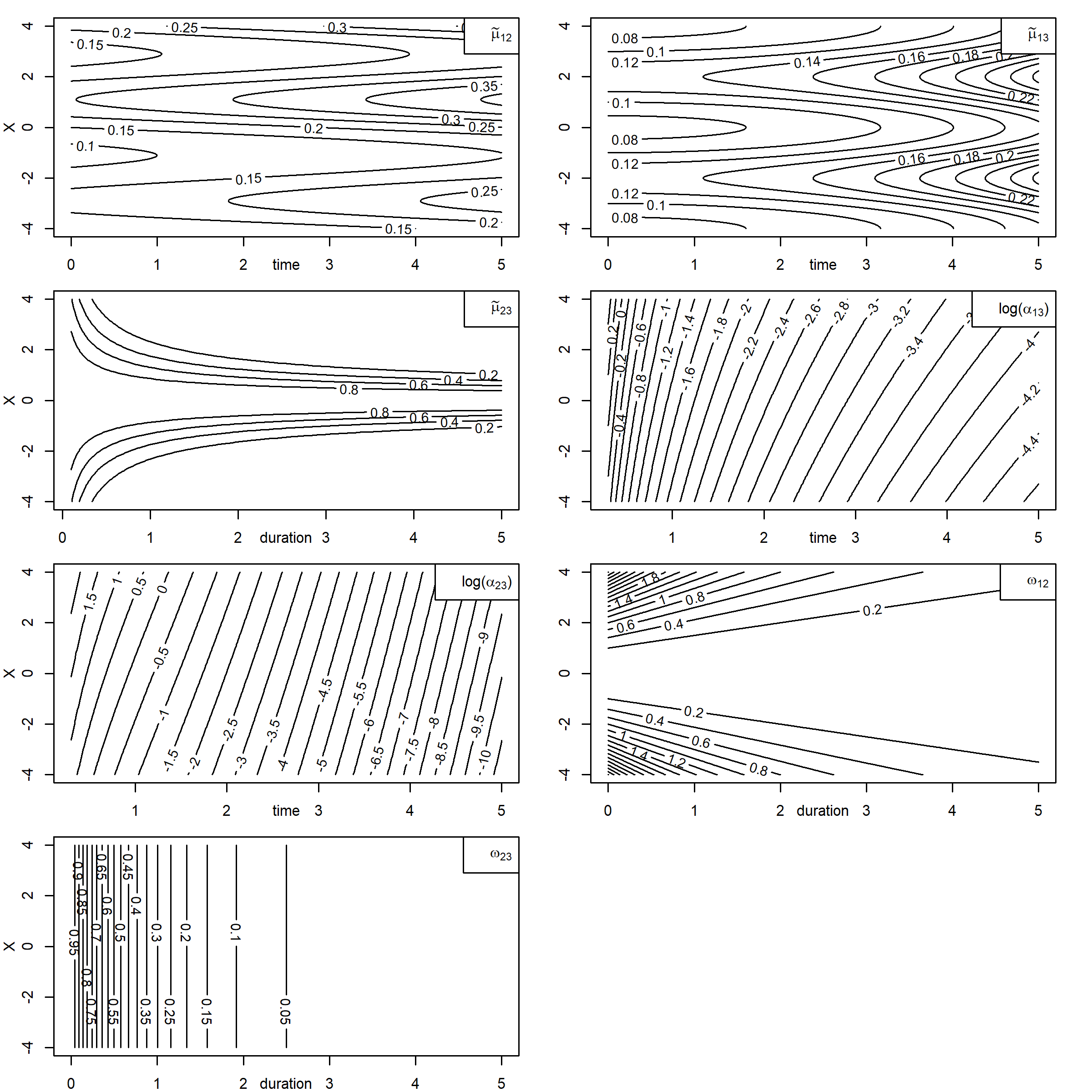}
    \caption{Contour plots of the hazards from the numerical study. The first axis is time or duration depending on whether the hazard depends on the time or the duration in the current state; none of the hazards depend on both. }
    \label{fig:Contour}
\end{figure}

Let the marked point process generated by the rates $\tilde{\mu}_{jk}$ be denoted $(T^\dagger_m,Y^\dagger_m)_{m \geq 1}$. Deletion of unreported events from $(T^\dagger_m,Y^\dagger_m)_{m \geq 1}$ leads to $(\tilde{T}_m,\tilde{Y}_m)_{m \geq 1}$, while deleting events that will be rejected leads to $(T^\ast_m,Y^\ast_m)_{m \geq 1}$. The data is generated according to the following algorithm:
\begin{enumerate}[(i)]
    \item Simulate $X$, $V$, and $C$. 
    \item Simulate events in $(V,C]$ using the rates $\tilde{\mu}_{jk}$ and initial state 1.
    \item For jumps to state 3, simulate corresponding reporting delays with a distribution that depends on whether or not a transition to state 2 has occurred previously.
    \item For transitions from state 2 to state 3, simulate adjudication events in $(T^\dagger_2+U^\dagger_2,\eta]$ if this interval is nonempty using the rates $\omega_{jk}$ and initial state 1.
\end{enumerate}
Events that are reported after $\eta$ are flagged and deleted except for when implementing oracle methods. The adjudication outcome $\xi_{\langle \eta \rangle}$ is simulated from a Bernoulli distribution with success probability equal to its expectation given $\mathcal{H}_\eta$, and events with adjudication outcome zero are deleted when implementing oracle methods.
For simulation of the reporting delays, we note that if $U \sim \textnormal{Weibull}(\lambda,k,\beta)$, then $[1-\exp\{-(\lambda U)^k\}]^{\exp(X \beta)} = W$ in distribution, where $W \sim \textnormal{Unif}(0,1)$. Trivial calculations then show that $U=\lambda^{-1} \times [-\log \{ 1-W^{\exp(-X \beta)}\}]^{1/k}$ in distribution. Contour plots of the hazards $\Tilde{\mu}_{jk}$, reporting delay reverse time hazards on log-scale, and adjudication hazards are depicted in Figure~\ref{fig:Contour}.

To obtain an expression for $\mu^\ast_{23}$, we first let $p(X)$ be the probability that a newly reported jump becomes confirmed. We have 
\begin{align*}
    p(X) &= \left(1-\exp\left[-\int_0^\infty g_{1} \times \{X/(v+2)\}^2 \diff v \right] \right)\left[1- \exp\left\{-\int_0^\infty  \exp(g_{2} \times v) \diff v \right\} \right] \\
    &= \left\{1-\exp\left(-g_{1} \times X^2/2\right)\right\}\{1-\exp(1/g_{2})\}.
\end{align*}
Note that $\Tilde{\mu}_{23}(t; \Tilde{T}_1, X)$ only depends on $t$ and $\Tilde{T}_1$ through $t-\Tilde{T}_1$, which equals $D(t)$ due to the lack of reporting delays for the transition to state 2, so we write $\Tilde{\mu}_{23}(t; \Tilde{T}_1, X)=\Tilde{\mu}_{23}\{D(t); X\}$. Using that thinning $(T^\dagger_m,Y^\dagger_m)_{m \geq 1}$ with the adjudication outcomes leads to $(T^\ast_m,Y^\ast_m)_{m \geq 1}$, similar calculations to the ones found in Section{\scshape~\ref{supplement:Distribution}} of the Supplementary material yield
\begin{align*}
    \mathbb{P}(T^\ast_2 \geq t \mid T^\ast_1, Y^\ast_1 = 2) &= 1-p(X)\left(1-\exp\left [ \frac{1-\exp\{\theta_7 D(t) X^2\}}{\theta_7 X^2} \right ] \right)
\end{align*}
and hence by Jacod's formula for the intensity, see Proposition~(3.1) of~\citet{Jacod:1975},
\begin{align*}
    \mu_{23}^\ast\{t,D(t);X,\boldsymbol{\theta}\} &= \frac{p(X) \exp\{ [1-\exp\{\theta_7 D(t) X^2\}]/(\theta_7 X^2)\}\Tilde{\mu}_{23}\{D(t);X\}}{1-p(X)(1-\exp\{[1-\exp\{\theta_7 D(t) X^2\}]/(\theta_7 X^2) \})}.
\end{align*}
In order to speed up the computation time when calculating $\hat{\boldsymbol\theta}_n$, we use that 
\begin{align*}
\int \mu^\ast_{1\bigcdot}(s; \emptyset, \boldsymbol\theta) \, \mathrm{d}s &= \left\{\pi/(4\theta_5)\right\}^{1/2} \exp\left\{\theta_4+\theta_6\times\cos(0.5 \pi X)\right\} \times \textnormal{erfi}( \theta_5^{1/2} \times s) \\
& \qquad +\theta_2^{-1} \times\exp\{\theta_1+\theta_2 \times X+\theta_3 \times \sin(0.5 \pi X)\} \times \exp(\theta_2 \times s)
\end{align*}
to calculate $Y_1(t) P^\ast \{ t; \mathbf{J}(t-), \boldsymbol\theta \}$ efficiently, where $\textnormal{erfi}$ is the imaginary error function. This function is implemented in many software packages, but we instead use the following 15'th order Taylor expansion around 0 to further speed up calculations: $\textnormal{erfi}(z) = \pi^{-1/2} \times (2z+\frac{2}{3}z^3+\frac{1}{5}z^5+\frac{1}{21}z^7+\frac{1}{108}z^9+\frac{1}{660}z^{11}+\frac{1}{4680}z^{13}+\frac{1}{37800}z^{15})+O(z^{17})$. The approximation error for $Y_1(t) P^\ast \{ t; \mathbf{J}(t-), \boldsymbol\theta \}$ was found to be of the order $10^{-7}$ for parameter choices near the true values, hence negligible. We further note that simple calculations gave analytical expressions for the exposure-terms of the log-likelihoods for $\mathbf{g}$ and $\mathbf{f}$ as well as for the Poisson approximation to the transition from state 1 to state 2.

\subsection{Results}

\begin{figure}[ht!]
    \centering
    \includegraphics[scale=0.725]{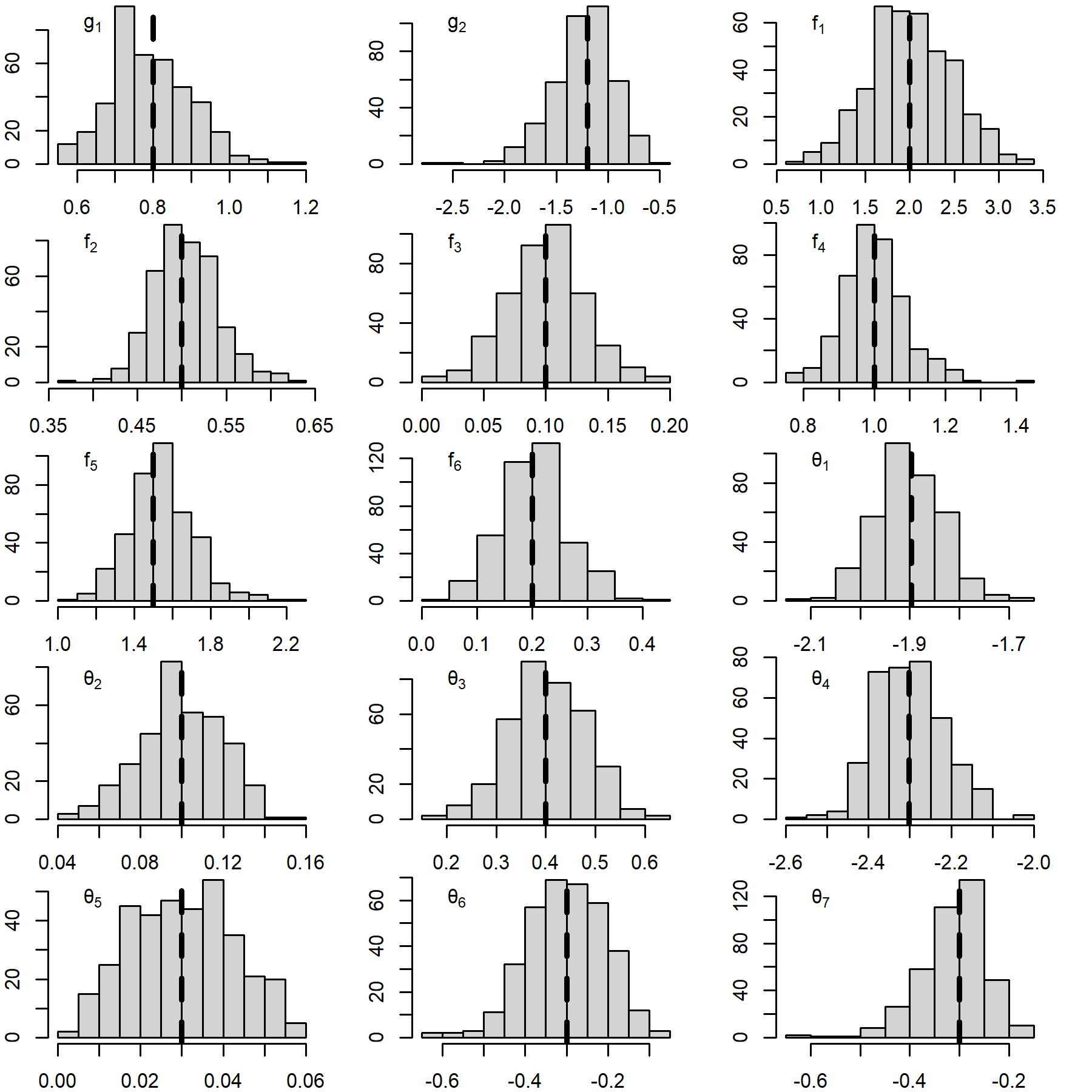}
    \caption{Histograms of $\hat{\mathbf{g}}_n$, $\hat{\mathbf{f}}_n$, and $\hat{\boldsymbol\theta}_n$ based on 400 samples of size $n=1500$ with the true values indicated by dashed lines.}
    \label{fig:histogramEstimators}
\end{figure}

In Table~\ref{table:gf} we report the bias, the SD, and the RMSE of the parameter estimators of $\mathbf{g}$ and $\mathbf{f}$.  In Table~\ref{table:thetarmse} we report the bias and the RMSE of the estimators of $\boldsymbol\theta$. The findings for SD and RMSE are highly comparable. Histograms of the obtained estimators are shown in Figure~\ref{fig:histogramEstimators} and generally show approximate Gaussianity with varying degrees of skewness. The bootstrap results are shown in Table~\ref{table:CI}, and are based on $k=50,100,200,300, 400$ estimates of $\theta_7$ and $1000$ bootstrap resamples. We report the coverage of the confidence bands for several values of the number of bootstrap runs $k$ and the confidence level $1-\alpha$.

\begin{table} 
\centering
\caption{Bias, empirical standard deviation (SD), and root mean squared error (RMSE) of the estimators $\hat{\mathbf{g}}_n$ and $\hat{\mathbf{f}}_n$ based on 400 simulations of size $n=1500$.}
\scalebox{0.9}{
\begin{tabular}{lccc}
\hline
Parameter &  Bias   & SD & RMSE \\[3pt]
\hline
$g_1=0.8$ & -.008  & .104 & .105   \\
$g_2=-1.2$ & -.035  & .294 & .296  \\
$f_1=2$ & .010  & .465 & .465  \\
$f_2=0.5$ & .004  & .037 & .037  \\
$f_3=0.1$  & -.000  & .032 & .032  \\
$f_4=1$  & -.001  & .089 & .089  \\
$f_5=1.5$  & .040  & .170 & .174  \\
$f_6=0.2$  & .004  & .060 & .060 \\
\hline
\end{tabular}}
\label{table:gf}
\end{table}

\begin{table} 
\centering
\caption{Bias and root mean squared error (RMSE) of the estimator $\hat{\boldsymbol\theta}_n$ based on 400 simulations of size $n=1500$.}
\scalebox{0.85}{
\begin{tabular}{lccccccccccccccc}
\hline
 &  \multicolumn{2}{c}{Proposed method}   &  \multicolumn{2}{c}{Oracle}  &  \multicolumn{2}{c}{Poisson approx.}  & \multicolumn{2}{c}{Naive 1} & \multicolumn{2}{c}{Naive 2}  \\
 \cmidrule(lr){2-3} \cmidrule(lr){4-5} \cmidrule(lr){6-7} \cmidrule(lr){8-9} \cmidrule(lr){10-11}
Parameter       &  Bias    & RMSE &  Bias    & RMSE  & Bias      & RMSE   & Bias    & RMSE   & Bias     & RMSE \\[3pt]
\hline
$\theta_1=\log 0.15$ & -.004    & .067 &  -.008   & .033  & -.010     & .067   & -.011   & .067   & -.010    &  .067  \\
$\theta_2 = 0.1$ & -.000    & .020 & -.001    & .020  & -.006     & .021   & -.006   & .021   & -.006    & .021 \\
$\theta_3 = 0.4$ & .003     & .078 & .003     & .078  & -.002     & .078   & -.002   & .078   & -.000    & .079 \\
$\theta_4 = \log 0.1$ & .003     & .083 & .001     & .083  & .012      & .092   & -.041   & .087   & -.051    & .096 \\
$\theta_5  =0.03$ & .000     & .012 & -.000    & .013  & -.006     & .017   & -.018   & .022   & -.015    & .021 \\
$\theta_6 = -0.3$ & -.000    & .094 & -.001    & .088  & .007      & .094   & -.009   & .087   & -.007    & .090 \\
$\theta_7 = -0.3$ & -.011    & .067 &  -.011   & .055  & -.012     & .067   & .157    & .158   & .148     & .163 \\
\hline
\end{tabular}}
\label{table:thetarmse}
\end{table}

\begin{table}
\centering
\caption{Coverage (\%) of percentile bootstrap confidence bands for $\theta_7$ with $k=50,100,200,300,400$ as well as $1-\alpha = 0.90,0.95,0.99.$}
\begin{tabular}{lcccccc}
\hline
$k$ &  \multicolumn{3}{c}{Coverage (\%)}    \\
\hline
    & 90  & 95 & 99   \\[3pt]    
\hline
50  & 91.8  & 93.9 & 95.9  \\ 
100  & 87.9  & 91.9 & 97.0  \\ 
200  & 87.9  & 93.0 & 98.0  \\ 
300  & 90.0  & 95.0 & 98.3  \\
400  & 89.5  & 95.0 &  98.7  \\
\hline
\end{tabular}
\label{table:CI}
\end{table}

\subsection{Empirical comparison with logistic regression approach}

In addition to the theoretical discussion of why~\citet{Cook:Kosorok:2004} would not be suitable for the application we have in mind, see Section~\ref{subsec:literature}, we also perform an empirical investigation in the same setting as the numerical study. As the Poisson approximation performed well in this setting, we also employ this approximation here to accommodate the effect of reporting delays. For adjudications, we like~\citet{Cook:Kosorok:2004} estimate the adjudication probabilities via logistic regression based on complete case data. Note that the adjudication hazards in the numerical study are always strictly positive (except on the null set where $X=0$). One interpretation of complete case data would therefore be those adjudications that are confirmed within the observation window, as these are the only adjudication outcomes that are known for certain based on the observed data. This would lead the estimated adjudication probabilities to be one for all of the events, which is clearly biased as it corresponds to using all unrefuted events. We instead treat all events that are not confirmed by time $\eta$ as rejected and regress the adjudication outcomes on $X$, $\tilde{T}_2$, and $\tilde{U}_2$. The only transitions affected by reporting delays are those from state $2$ to state $3$, so we only estimate $f_4$,$f_5$,$f_6$, and $\theta_7$. The results are given in Table~\ref{table:cookkosorok}. With the exception of the $f_5$ estimate, this approach leads a large increase in bias for the parameter estimates, with the performance for $\theta_7$ being substantially worse than both Naive 1 and Naive 2.

\begin{table}[ht!]
\centering
\caption{Bias and empirical standard deviation (SD) for estimators of $f_4$, $f_5$, $f_6$ and $\theta_7$ based on 400 simulations of size $n=1500$.}
\scalebox{0.9}{
\begin{tabular}{lcc}
\hline
Parameter &  Bias   & SD \\[3pt]
\hline
$f_4=1$  & .151  & .156   \\
$f_5=1.5$  & -.013  & .111   \\
$f_6=0.2$  & .015  & .052  \\
$\theta_7=-0.3$  & -.425  & .123  \\
\hline
\end{tabular}}
\label{table:cookkosorok}
\end{table}

\section{Additional Details for Data Application} \label{supplement:Application}

In Table~\ref{table:dataf} and Table~\ref{table:datag}, we present the parameter estimates and percentile bootstrap confidence intervals for $\mathbf{f}$ and $\mathbf{g}$, respectively. The confidence intervals are computed using 400 bootstrap resamples. The variable `rejected before' is an indicator of the event that adjudication state 2 has been entered previously.  This variable is trivial in adjudication state 2 and hence not included in that regression. There are no deaths observed from adjudication state 2 for disablements, and we hence estimate no parameters for this transition, letting the hazard be identically zero. Similarly, only one death is observed from adjudication state 1 for reactivations, causing the Poisson regression algorithm to diverge, and we hence also let this hazard be identically zero. One could of course have changed the regression model to a simple intercept model and obtained a small positive hazard, but its effect on the adjudication probabilities would be negligible and we hence forego this modification.

\begin{table}[ht!]
\centering
\scalebox{0.93}{
{%
\Rotatebox{90}{%
\begin{tabular}{llcccccccccc}
 & &  \multicolumn{2}{c}{$\omega_{12}$} &  \multicolumn{2}{c}{$\omega_{13}$} &  \multicolumn{2}{c}{$\omega_{14}$} &  \multicolumn{2}{c}{$\omega_{21}$} &  \multicolumn{2}{c}{$\omega_{24}$}   \\
 \hline
Parameter & Type &  Est.   & CI &  Est.   & CI &  Est.   & CI &  Est.   & CI &  Est.   & CI  \\[3pt]
\hline
Age & Reactivation & .005  & (-0.010, 0.018)  &  .015  & (0.007, 0.026)   &  -  &  -  &  .000  &  (-0.011, 0.010)   &  .023  & (-0.059, 0.105)  \\
Male & Reactivation & .973  & (0.295, 1.63)  &  1.29  & (0.752, 1.76)  &  -  & -   &  -.208  &  (-0.703, 0.295)   &  -5.94  &  (-10.45, -2.87)  \\
Female & Reactivation & 1.18  & (0.524, 1.84)  &  .771  & (0.159, 1.31)  &  -  & -   &  -.053  &  (-0.793, 0.626)   &  -21.10  & (-26.24, -18.17)  \\
Disability duration & Reactivation & -.065  & (-0.247, 0.153)  &  -.198  & (-0.380, 0.005)   &  -  &  -   &  .206  &  (0.101, 0.325)  &  -.464  &  (-1.86, 0.136)  \\
Reactivation duration & Reactivation & -.042  & (-0.273, 0.205)  &  -.490  & (-0.803, -0.200)   &  -  & -   &  -1.25  &  (-1.49, -1.05)  &  .148  &  (-1.29, 1.66)  \\
Age & Disability & -.005  & (-0.015, 0.005) & .001 & (-0.003, 0.004)   &  .013  &  (-0.054, 0.101)  & .003  &  (-0.008, 0.014)  &  -  & - \\
Male & Disability & -.800  & (-1.27, -0.310) & 1.06 &  (0.870, 1.25)  &  -4.85  &  (-10.04, 4.52)  & 1.16  &  (0.578, 1.74)  &  -  & - \\
Female & Disability & -.806  & (-1.30, -0.316) & 1.19 &  (0.988, 1.41)  &  -5.29  &  (-22.06, 3.17)  & 1.08  &  (0.353, 1.71)  &  -  & - \\
Disability duration & Disability & .077  & (-0.190, 0.315) & -.234 &  (-0.343, -0.129)  &  .169  &  (-1.45, 1.23)  & -.255  &  (-0.610, 0.063)  &  - & - \\
Report duration & Disability & -.250  & (-0.802, 0.280) & .470 & (0.307, 0.680)   &  -.532  &  (-211.85, 2.12)  & -.090  &  (-0.525, 0.413)  &  -  & - \\
Rejected before & Disability & .982  & (0.572, 1.38) & -.118 &  (-0.360, 0.072)  &  -14.68  &  (-16.70, -11.75)  & -  &  -  &  -  & - \\
\hline
\end{tabular}
}}}
\caption{Parameter estimates (Est.) and 95\% bootstrap percentile confidence interval (CI) for the parameter $\mathbf{g}$ using the proposed method and 400 bootstrap resamples.}
\label{table:datag}
\end{table}

\begin{table}[ht!]
\centering
\caption{Parameter estimates (Est.) and 95\% bootstrap percentile confidence interval (CI) for the parameter $\mathbf{f}$ using the proposed method and 400 bootstrap resamples.}
\begin{tabular}{llcc}
\hline
 Parameter &   Est.   & CI   \\[3pt]
 \hline
$\lambda$ & 2.23 & (0.964, 2.78)  \\
$k$ & 1.05 & (0.931, 2.00)     \\
Male & .356 & (-0.518, 0.682)    \\
Female & .113 & (-0.672, 0.783) \\
Age at disability & .001 & (-0.027, 0.006) \\
\hline
\end{tabular}
\label{table:dataf}
\end{table}

\FloatBarrier

\noindent To further illustrate the combined effect of the estimated hazards, we calculate the transition probabilities 
\begin{align*}
(t,j,a) \mapsto \mathbb{P}_{\hat{\boldsymbol\theta}_n}\{ I^\ast_j(t) \mid I^\ast_1(\eta)=1, \textnormal{Age}(\eta)=30,\textnormal{Gender}=a\}
\end{align*}
for $t \geq \eta$. They represent the probability that someone who at the end of the observation window is active, 30 years old, and with gender $a$, will be in state $j$ at a later time. The data does not permit a reasonable estimation of mortalities (death hazards), so we simply employ the mortalities from Table~\ref{table:deathHazards}, which are inspired by actuarial practice and, at least, admit a reasonable level. Like the disability and reactivation hazards, the covariates in the mortality hazards enter in a linear predictor with log link. To not overextrapolate the estimated calendar time effect for the disability hazard, we keep its effect fixed at its value at time $\eta$. The transition probabilities are calculated by solving the differential equation from Theorem 3.1 in~\citet{Buchardt.etal:2015} using an Euler scheme, see Section 3 in their paper for details on the algorithm. 

\begin{figure}[ht!]
    \centering
    \includegraphics[scale=0.55]{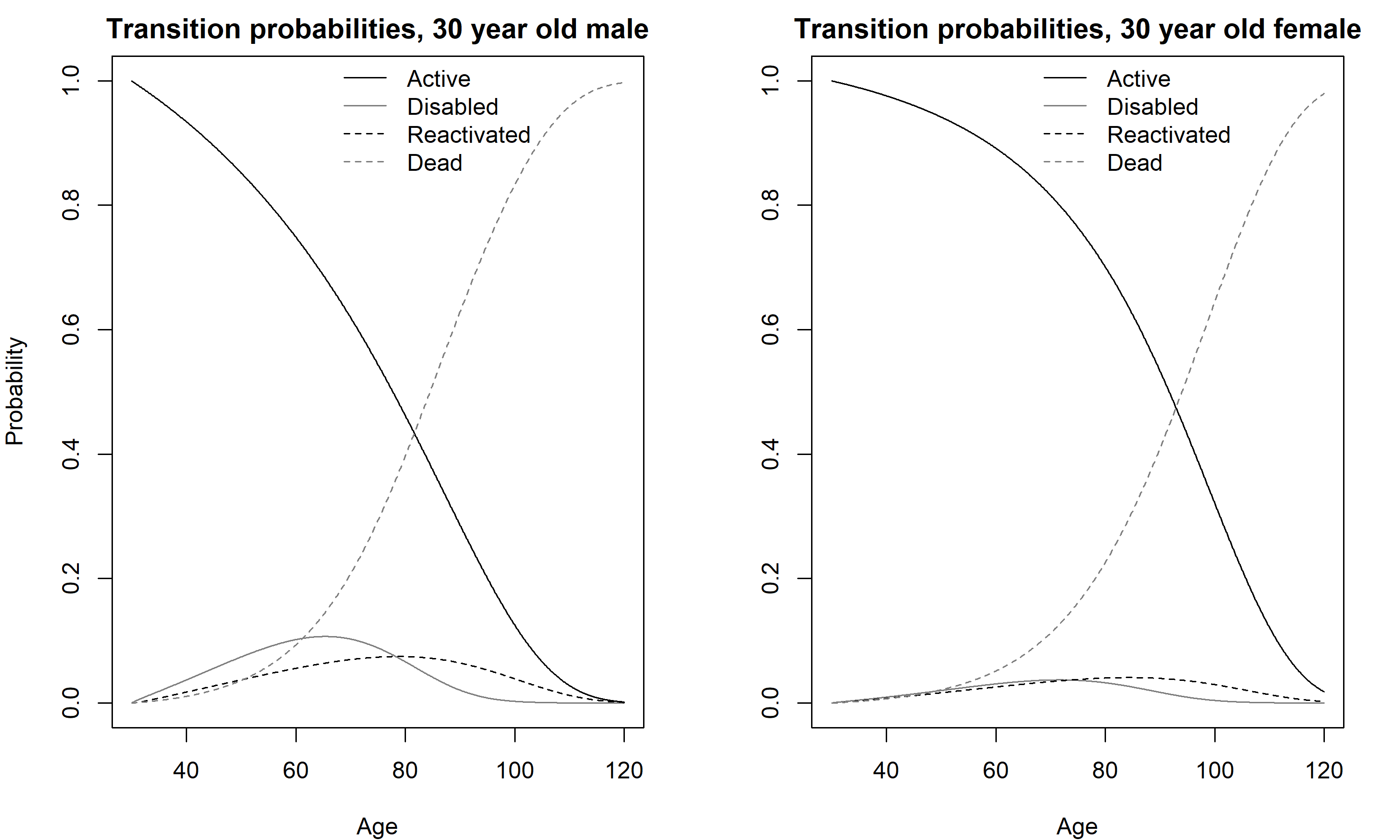}
    \caption{Transition probabilities from age 30 to age 120 for a male (left) and female (right) that is active and 30 years old at time $\eta$. }
    \label{fig:Probabilities}
\end{figure}

The results are depicted in Figure~\ref{fig:Probabilities}. As the estimated disability hazards would suggest, males are generally more likely to become disabled than females. The proportion of disabled people is larger than the proportion who have died for the first few years, but is overtaken around ages 50 to 60. The proportion of reactivated individuals peaks around the ages 80 to 90 and decreases afterwards as the death probability begins to dominate.

\begin{table}[ht!]
\centering
{%
\caption{Mortalities inspired by actuarial practice.}
\begin{tabular}{lcccc}
\hline
 &  Age  & Male & Female & $\min\{\textnormal{Duration},5\}$ \\[3pt]
\hline
$\mu_{ad}$, $\mu_{rd}$ & 0.07  & -9.50 & -9.80 & --     \\
$\mu_{id}$ & 0.07 & -6.40 & -6.80 & -0.25 \\
\hline
\end{tabular}}
\label{table:deathHazards}
\end{table}

\FloatBarrier